\documentclass[review]{elsarticle}

\usepackage{lineno,hyperref}
\modulolinenumbers[5]

\journal{Wave Motion}

\biboptions{sort&compress}

\usepackage{amsmath}
\usepackage{amsfonts}

\usepackage{cleveref}
\crefname{equation}{}{}
\crefname{figure}{Figure}{Figure}
\crefname{section}{Section}{Section}
\usepackage{bm}

\usepackage{color}
\usepackage{ulem}

\makeatletter
\let\MYcaption\@makecaption
\makeatother
\usepackage{subcaption}
\makeatletter
\let\@makecaption\MYcaption
\makeatother

\newcommand{\ione}{\mathrm{i}}
\newcommand{\e}{\mathrm{e}}
\newcommand{\dd}{\mathrm{d}}
\newcommand{\eps}{\varepsilon}
\renewcommand{\Re}{\mathrm{Re}\,}
\renewcommand{\Im}{\mathrm{Im}\,}

\begin{document}

\begin{frontmatter}

\title{A topology optimization of open acoustic waveguides based on a scattering matrix method}


\author[tokyo]{Kei Matsushima\corref{mycorrespondingauthor}}
\cortext[mycorrespondingauthor]{Corresponding author}

\author[keio]{Hiroshi Isakari}

\author[nagoya]{Toru Takahashi}

\author[nagoya]{Toshiro Matsumoto}

\address[tokyo]{The University of Tokyo, 2-11-16 Yayoi, Bunkyo-ku, Tokyo, Japan}
\address[keio]{Keio University, 3-14-1 Hiyoshi, Kohoku-ku, Yokohama, Kanagawa, Japan}
\address[nagoya]{Nagoya University, Furo-cho, Chikusa-ku, Nagoya, Aichi, Japan}

\begin{abstract}
This study presents a topology optimization scheme for realizing a bound state in the continuum along an open acoustic waveguide comprising a periodic array of elastic materials. First, we formulate the periodic problem as a system of linear algebraic equations using a scattering matrix associated with a single unit structure of the waveguide. The scattering matrix is numerically constructed using the boundary element method. Subsequently, we employ the Sakurai--Sugiura method to determine resonant frequencies and the Floquet wavenumbers by solving a nonlinear eigenvalue problem for the linear system. We design the shape and topology of the unit elastic material such that the periodic structure has a real resonant wavenumber at a given frequency by minimizing the imaginary part of the resonant wavenumber. The proposed topology optimization scheme is based on a level-set method with a novel topological derivative. We demonstrate a numerical example of the proposed topology optimization and show that it realizes a bound state in the continuum through some numerical experiments.
\end{abstract}


\end{frontmatter}


\section{Introduction}
Recently, bound states in the continuum (BICs) have been enthusiastically investigated in the fields of quantum mechanics, photonics, {acoustics, and water waves} \cite{hsu2016bound}. BICs were originally proposed by von Neumann and Wigner \cite{vonneumann1929uber} in a quantum system and then experimentally found in some classical systems \cite{plotnik2011experimental,lee2012observation,hsu2013observation}. Historically, it has been said that bound states can exist only outside the radiation continuum, where a near-field state cannot be coupled to any radiation (scattering) channel, resulting in a perfectly confined state in the vicinity of a structure. Challenging this conventional wisdom, BICs may exist within the continuum in some geometrical systems, e.g., waveguides \cite{linton2007embedded,dreisow2009adiabatic} and photonic/phononic crystal slabs (diffraction gratings) \cite{lee2012observation,hsu2013observation,yang2014analytical,bulgakov2014bloch}. BICs are of theoretical interest and practical importance due to their potential to realize high-Q resonance, which is an essential property of lasers \cite{kodigala2017lasing}, filters \cite{foley2014symmetryprotected}, and sensors \cite{yanik2011seeing} for next generations.

Because resonance properties are sensitive to the material and geometrical configurations of a structure, some inverse-design approaches, such as parameter tuning and shape/topology optimization \cite{sokolowski1992introduction,bendsoe2013topology}, may be necessary to realize high-Q resonance originating from BICs in practical applications {\cite{hsu2016bound}}. Such optimization techniques have been recently used to design some photonic and phononic structures with maximized bandgaps \cite{dobson1999maximizing}. Though BICs are formulated similar to photonic/phononic bandgaps, we encounter some numerical difficulties when applying {optimization-based design methods to manipulate BICs because they often rely on finite element methods}.

BICs are formulated as resonant states satisfying Maxwell's or Helmholtz' equations in open systems and characterized as an eigenmode of a nonlinear eigenvalue problem to find a frequency and wavenumber that allow a nonzero state without any incident field. For some simple geometries, we can employ a semianalytical technique, such as cylindrical or spherical wave expansions \cite{bulgakov2014bloch}, to effectively compute resonant states. {For example, Evans and Porter showed numerical evidence that a circular inclusion in a planar waveguide supports BICs \cite{evans1998trapped}. Later, rectangular scatterers are considered by Porter and Evans \cite{porter2005embedded}. The recent work by Bennetts and Peter investigated arrays of circular scatterers in water based on a transfer operator with cylindrical functions \cite{bennetts2022rayleigh}. See the review by Linton and McIver for more details \cite{linton2007embedded}.} However, more general configurations require a discretization-based method, such as the finite element method and boundary element method (BEM). These methods should be carefully implemented to avoid neglecting the radiation effect of BICs because resonant states in open systems are not bounded or even diverging in space \cite{hu2009understanding}. This is not the case for the photonic bandgap computation because the underlying eigenvalue problem is defined in a bounded domain.

After the numerical analysis of BICs, we need to calculate the sensitivity of their eigenvalues for structural optimization, which is called design sensitivity. Because the corresponding boundary value problem is defined in an open space, we need to truncate the unbounded domain and to carefully deal with the truncated boundary through some special treatment, such as the Dirichlet-to-Neumann map to evaluate the variation in the eigenvalue with respect to geometrical perturbation \cite{ammari2020perturbation}. This would incur additional computational cost (especially when the BEM is used) because the variation involves a volume integral of the resonant state over the truncated domain. To the authors' best knowledge, no prior works have found design sensitivity for BICs.

Therefore, this study proposes a topology optimization scheme for designing open acoustic waveguides exhibiting BICs at desired frequencies in two dimensions to overcome the above difficulties. First, we propose a scattering matrix-based approach to compute BICs. The basic idea is the same as the one proposed in \cite{bulgakov2014bloch}, which calculates BICs along a periodic array of circular rods. For the topology optimization, we incorporate the BEM into the scattering matrix method to deal with more geometrically complex structures than the circular rods \cite{gimbutas2013fast}. Further, we formulate the topological derivative \cite{sokolowski1999topological} of a resonant wavenumber based on the scattering matrix method. This formulation does not require any volume integration of resonant states; thus, it saves considerable computational costs. Subsequently, we incorporate the topological derivative into a level set-based topology optimization algorithm \cite{isakari2017topology}. Finally, we perform the topology optimization and numerically demonstrate that the optimized structure forms a BIC. 

\section{Scattering by multiple and periodic obstacles}
In this section, we first formulate wave scattering by multiple and periodic obstacles in two dimensions using the scattering matrix method. Further, we describe how to compute the eigenvalues of the systems.

\subsection{Scattering through a single obstacle}
\begin{figure}[h]
  \begin{center}
  \includegraphics[scale=0.4]{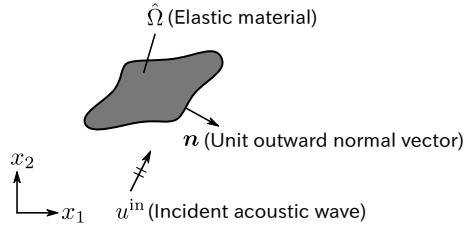}
   \caption{Scattering through a single scatterer $\hat{\Omega}$ placed in the two-dimensional space $\mathbb{R}^2$.}
   \label{fig:single}
\end{center}
\end{figure}
As shown in \cref{fig:single}, we first consider a scattering problem where a single elastic material $\hat{\Omega}$ is placed in the free space $\mathbb{R}^2$. Throughout this paper, we neglect the shear modulus and formulate the time-harmonic scattering problem using the following transmission problem:
\begin{align}
 \nabla^2 u(\bm{x}) + \frac{\omega^2}{c^2} u(\bm{x}) = 0 &\quad x\in\Omega:=\mathbb{R}^2\setminus\overline{\hat{\Omega}}, \label{eq:bvp_1}
 \\
 \nabla^2 u(\bm{x}) + \frac{\omega^2}{\hat{c}^2} u(\bm{x}) = 0 &\quad x\in\hat{\Omega}, \label{eq:bvp_2}
 \\
 u(\bm{x}):=u|_+(\bm{x}) = u|_-(\bm{x}) &\quad x\in\partial\hat{\Omega}, \label{eq:bvp_3}
 \\
 q(\bm{x}):=\frac{1}{\rho}\frac{\partial u}{\partial n}\Big|_+(\bm{x}) = \frac{1}{\hat{\rho}}\frac{\partial u}{\partial n}\Big|_-(\bm{x}) &\quad x\in\partial\hat{\Omega}, \label{eq:bvp_4}
 \\
 \left| \frac{\partial}{\partial r}(u-u^\mathrm{in}) - \ione k(u-u^\mathrm{in}) \right| = O(r^{-1/2}) \ \mathrm{as} &\quad r=|\bm{x}|\to\infty, \label{eq:bvp_5}
\end{align}
where $u$ denotes the sound pressure, $k=\omega/c$ the wavenumber, and $u^\mathrm{in}$ the corresponding incident wave. {The overline denotes the closure of a domain.} {A vector quantity $\bm{p}$ associated with the Cartesian coordinate system $(x_1,x_2)$ is denoted by a bold symbol, and its components are expressed by $p_i$ ($i=1,2$).} In addition, $\frac{\partial }{\partial n}=\bm{n} \cdot \nabla$ denotes the normal derivative with the unit outward normal vector $\bm{n}$ to $\hat{\Omega}$. \Cref{eq:bvp_5} represents the Sommerfeld radiation condition. Throughout the paper, the time dependence of the time-harmonic fields is chosen as $\e^{-\ione\omega t}$ with the angular frequency $\omega$. In the exterior medium $\Omega$ and scatterer $\hat{\Omega}$, the phase velocities $c$ and $\hat{c}$ are given by their mass densities $\rho$ and $\hat{\rho}$ and bulk moduli $\kappa$ and $\hat{\kappa}$ as $c=\sqrt{\kappa/\rho}$ and $\hat{c}=\sqrt{\hat{\kappa}/\hat{\rho}}$, respectively. The symbol $|_+$ (resp. $|_-$) denotes the trace from $\Omega$ (resp. $\hat{\Omega}$) to the boundary.

\subsection{Scattering matrix}

\subsubsection{Definition}
Using Graf's addition theorem, which is given by
\begin{align}
 &\quad H^{(1)}_n(k|\bm{x}-\bm{y}|)\e^{\ione n\theta (\bm{x}-\bm{y})} 
 \notag\\
 &= 
 \begin{cases}
  \displaystyle\sum_{m=-\infty}^\infty H^{(1)}_{m}(k|\bm{x}|)J_{m-n}(k|\bm{y}|)\e^{\ione m\theta(\bm{x})}\e^{-\ione(m-n)\theta(\bm{y})} & (|\bm{x}|>|\bm{y}|)
  \\
  \displaystyle\sum_{m=-\infty}^\infty J_{m}(k|\bm{x}|)H^{(1)}_{m-n}(k|\bm{y}|)\e^{\ione m\theta(\bm{x})}\e^{-\ione(m-n)\theta(\bm{y})} & (|\bm{x}|<|\bm{y}|)
 \end{cases}
 \label{eq:graf}
\end{align}
for the Hankel function $H^{(1)}_n$ of the first kind and Bessel function $J_n$ of order $n$ with $\theta(\bm{x})=\tan^{-1} (x_2/x_1)$, we have the multipole expansion of the fundamental solution $G(\bm{x},\bm{y})=\frac{\ione}{4}H^{(1)}_0(k|\bm{x}-\bm{y}|)$ as
\begin{align}
 G(\bm{x},\bm{y}) = \frac{\ione}{4}\sum_{n=-\infty}^\infty H^{(1)}_{n}(k|\bm{x}-\bm{x}_0|)J_n(k|\bm{y}-\bm{x}_0|)\e^{\ione n\theta(\bm{x}-\bm{x}_0)}\e^{-\ione n\theta(\bm{y}-\bm{x}_0)}
\end{align}
for any $\bm{x}_0\in\mathbb{R}^2$ and $(\bm{x},\bm{y})$ satisfying $|\bm{x}-\bm{x}_0|>|\bm{y}-\bm{x}_0|$. Substituting this series into the following representation formula
\begin{align}
 u(\bm{x}) = u^\mathrm{in}(\bm{x}) - \rho \int_{\partial \hat{\Omega}} G(\bm{x},\bm{y})q(\bm{y})\dd\Gamma_y + \int_{\partial \hat{\Omega}} \frac{\partial G}{\partial n_y}(\bm{x},\bm{y})u(\bm{y})\dd\Gamma_y \quad \bm{x}\in\Omega,
\end{align}
we obtain that the solution $u$ can also be written in terms of the cylindrical functions as follows:
\begin{align}
 u(\bm{x}) = u^\mathrm{in}(\bm{x}) + \sum_{n=-\infty}^\infty B_n O_n(\bm{x}-\bm{x}_0) \quad \bm{x}\in\Omega\setminus\overline{D}, \label{eq:outgoing}
 \\
 B_n = \frac{\ione(-1)^n}{4} \langle I_{-n}, u \rangle_{\partial\hat{\Omega}}, \label{eq:B_integ}
 \\
 I_n(\bm{x}) = J_n(k|\bm{x}|)\e^{\ione n\theta(\bm{x})},
 \\
 O_n(\bm{x}) = H^{(1)}_n(k|\bm{x}|)\e^{\ione n\theta(\bm{x})},
\end{align}
where $\langle \cdot,\cdot \rangle$ denotes the bilinear form defined by the following boundary integral:
\begin{align}
 \langle v, u \rangle_{\partial\hat{\Omega}} = \int_{\partial\hat{\Omega}}\left(\frac{\partial v}{\partial n}\bigg|_+ u|_+ - v|_+\frac{\partial u}{\partial n} \bigg|_+\right)\dd\Gamma, 
\end{align}
and $D=\{\bm{x}\in\mathbb{R}^2 \mid |\bm{x}-\bm{x}_0|< \max_{\bm{y}\in \partial\hat{\Omega}} |\bm{y}-\bm{x}_0|  \}$ denotes the minimum enclosing disk of $\hat{\Omega}$ centered at $\bm{x}_0$. The representation \cref{eq:outgoing} implies that the coefficient vector $B=(B_n)_{n\in\mathbb{Z}}$ completely describes the field $u$ in the exterior of the disk $D$. 

We suppose that the incident wave $u^\mathrm{in}$ can also be expanded into the cylindrical functions as follows:
\begin{align}
 u^\mathrm{in}(\bm{x}) = \sum_{n=-\infty}^\infty A_n I_n(\bm{x}-\bm{x}_0) \quad \bm{x}\in\Omega\setminus\overline{D} \label{eq:incoming}
\end{align}
with complex coefficients $A_n\in\mathbb{C}$. For example, when $u^\mathrm{in}$ is the plane wave $\e^{\ione k \bm{p}\cdot (\bm{x}-\bm{x}_0)}$ propagating along a unit vector $\bm{p}\in\mathbb{R}^2$, we have $A_n=(p_2+\ione p_1)^n$. From the linearity of the scattering problem \cref{eq:bvp_1,eq:bvp_2,eq:bvp_3,eq:bvp_4,eq:bvp_5}, the relationship between the incident coefficient vector $A=(A_n)_{n\in\mathbb{Z}}$ and $B$ should also be linear, which yields the linear equation
\begin{align}
 B_n = \sum_{n^\prime=-\infty}^\infty S_{nn^\prime} A_{n^\prime}. \label{eq:S_def}
\end{align}
The matrix $S$ is called a {\it scattering matrix}. In what follows, such equation is simply denoted using the matrix-vector notation $B=SA$. {For a given vector $A$, we can compute the multiplication $B=SA$ as follows:
\begin{enumerate}
	\item Set the incident wave $u^\mathrm{in}(\bm{x}) = \sum_{n=-\infty}^\infty A_n I_n(\bm{x}-\bm{x}_0)$.
	\item Solve the scattering problem (1)--(5) for a given shape $\hat{\Omega}$ and compute $u$ and $q$ on $\partial \hat{\Omega}$.
	\item Compute the boundary integral \cref{eq:B_integ}.
\end{enumerate}
See \cite{gimbutas2013fast} for details.
}

\subsubsection{Boundary element method}
Once we have a scattering matrix for the scattering problem \cref{eq:bvp_1,eq:bvp_2,eq:bvp_3,eq:bvp_4,eq:bvp_5}, the scattered field is uniquely determined using \cref{eq:outgoing,eq:S_def}.  From the definition \cref{eq:S_def}, each component of the scattering matrix $S$ is given by $S_{nn^\prime}=\langle I_{-n}, u_{n^\prime} \rangle_{\partial\hat{\Omega}}$, where $u_{n^\prime}$ denotes the solution of the boundary value problem \cref{eq:bvp_1,eq:bvp_2,eq:bvp_3,eq:bvp_4,eq:bvp_5} for $u^\mathrm{in}=I_{n^\prime}$. Such a solution can be obtained by solving an appropriate boundary integral equation. In this study, we use the Burton--Miller-type boundary integral equation \cite{burton1971application}, which is given by
\begin{align}
 \begin{bmatrix}
  \frac{1}{2}\mathcal{I} - \mathcal{D} - \eta\mathcal{N} & \rho\left(\mathcal{S} +\eta\left(\frac{1}{2}\mathcal{I}+\mathcal{D}^*\right) \right)
  \\
  \frac{1}{2}\mathcal{I} + \hat{\mathcal{D}}             &  -\hat{\rho}\hat{\mathcal{S}}
 \end{bmatrix}
 \begin{pmatrix}
  u_{n^\prime} \\ q_{n^\prime}
 \end{pmatrix}
 =
 \begin{pmatrix}
  I_{n^\prime} + \eta \frac{\partial I_{n^\prime}}{\partial n}
  \\
  0
 \end{pmatrix}
 ,
\end{align}
for $u_{n^\prime}$ and $q_{n^\prime}=\frac{1}{\rho}\frac{\partial u_{n^\prime}}{\partial n} |_+$ with coupling parameter $\eta\in\mathbb{C}$, where the integral operators $\mathcal{S}$, $\mathcal{D}$, $\mathcal{D}^*$, and $\mathcal{N}$ are respectively defined as
\begin{align}
 (\mathcal{S}\phi)(\bm{x}) &= \int_{{\partial \hat{\Omega}}} G(\bm{x},\bm{y})\phi(\bm{y})\dd\Gamma_y,
 \\
 (\mathcal{D}\phi)(\bm{x}) &= \int_{{\partial \hat{\Omega}}} \frac{\partial G}{\partial n_y}(\bm{x},\bm{y})\phi(\bm{y})\dd\Gamma_y,
 \\
 (\mathcal{D}^*\phi)(\bm{x}) &= \int_{{\partial \hat{\Omega}}} \frac{\partial G}{\partial n_x}(\bm{x},\bm{y})\phi(\bm{y})\dd\Gamma_y,
 \\
 (\mathcal{N}\phi)(\bm{x}) &= \mathrm{p.f.}\,\int_{{\partial \hat{\Omega}}} \frac{\partial^2 G}{\partial n_x \partial n_y}(\bm{x},\bm{y})\phi(\bm{y})\dd\Gamma_y,
\end{align}
and ``$\mathrm{p.f.}$'' represents the finite part of the divergent integral. Moreover, the other operators $\hat{\mathcal{S}}$ and $\hat{\mathcal{D}}$ are obtained by replacing $c$ in $\mathcal{S}$ and $\mathcal{D}$ with $\hat{c}$, respectively. {The parameter $\eta\in\mathbb{C}$ is introduced to avoid fictitious eigenvalues, at which the the boundary integral equation becomes ill-posed \cite{burton1971application}. Although $\eta$ is arbitrary as long as $\Im[\eta]<0$, the formula $\eta=-\ione/k$ is known to be the best choice in terms of the condition number of a discretized system \cite{zheng2015burton}}.

\subsection{Scattering through a finite number of obstacles}\label{ss:multi}
Next, we describe a scattering-matrix formalism to solve multiple scattering problems. {See \cite{abramowitz1965handbook,martin2006multiple} for details.}

\begin{figure}[h]
  \begin{center}
  \includegraphics[scale=0.4]{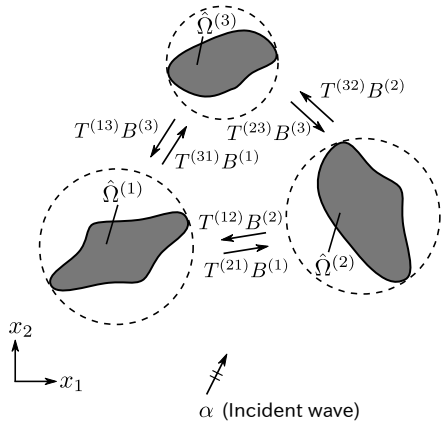}
   \caption{Scattering through multiple scatterers $\hat{\Omega}^{(i)}$ placed in the two-dimensional space $\mathbb{R}^2$}
   \label{fig:multi}
\end{center}
\end{figure}
As shown in \cref{fig:multi}, we consider the scattering through $N$ scatterers $\hat{\Omega}^{(i)}$ ($i=1,\ldots,N$). The shapes and materials of the scatterers are not necessarily identical. Let $S^{(i)}$ denote the scattering matrix associated with the scatterer $\hat{\Omega}^{(i)}$. The only assumption here is that any minimum disk $D^{(i)}$ enclosing $\hat{\Omega}^{(i)}$, whose center is denoted by $x_0^{(i)}$, does not overlap with each other (well-separated condition).

Under this assumption, we can write the total field as follows:
\begin{align}
 u(\bm{x}) = u^\mathrm{in}(\bm{x}) + \sum_{i=1}^N \sum_{n=-\infty}^\infty B^{(i)}_n O_n(\bm{x}-\bm{x}_0^{(i)}) \quad \bm{x}\in \mathbb{R}^2\setminus\overline{\cup_{i=1}^N D^{(i)}}
\end{align}
for outgoing multipole coefficients $B^{(i)}=(B_n^{(i)})_{n\in\mathbb{Z}}$ associated with $\hat{\Omega}^{(i)}$. Our task is to describe relations among $B^{(1)},\ldots,B^{(N)}$ using the scattering matrices $S^{(i)}$. If $\bm{x}$ is located around $D^{(j)}$, i.e., $\bm{x}\in\mathbb{R}^2\setminus\overline{\cup_{i=1}^N D^{(i)}}$ and $\forall i,\ |\bm{x}-\bm{x}_0^{(j)}|<|\bm{x}_0^{(i)}-\bm{x}_0^{(j)}|$ hold, then we can use the formula \cref{eq:graf} to obtain
\begin{align}
 & \sum_{n=-\infty}^\infty B^{(i)}_n O_n(\bm{x}-\bm{x}_0^{(i)})
 \notag\\
 &= \sum_{n=-\infty}^\infty B^{(i)}_n H^{(1)}_n(k|\bm{x}-\bm{x}_0^{(j)} - (\bm{x}_0^{(i)}-\bm{x}_0^{(j)})|) \e^{\ione n\theta(\bm{x}-\bm{x}_0^{(j)} - (\bm{x}_0^{(i)}-\bm{x}_0^{(j)}))}
 \notag\\
 &= \sum_{n=-\infty}^\infty \left( \sum_{m=-\infty}^\infty H^{(1)}_{m-n}(k|\bm{x}_0^{(j)}-\bm{x}_0^{(i)}|)\e^{\ione (m-n)\theta(\bm{x}_0^{(j)}-\bm{x}_0^{(i)})} B^{(i)}_n\right) J_n(k|\bm{x}-\bm{x}_0^{(j)}|)\e^{\ione n\theta(\bm{x}-\bm{x}_0^{(j)})}
 \notag\\
 &= \sum_{n=-\infty}^\infty \left(T^{(ji)} B^{(i)}\right)_n I_n(\bm{x}-\bm{x}_0^{(j)}),
 \label{eq:translation}
\end{align}
where $T^{(ji)}_{nm}=O_{m-n}(\bm{x}_0^{(j)}-\bm{x}_0^{(i)})$ denotes a translation matrix from $i$th to $j$th scatterer. The translation formula \eqref{eq:translation} indicates that the outgoing wave $B^{(i)}$ from $\hat{\Omega}^{(i)}$ turns into the incoming wave $T^{(ji)}B^{(i)}$ around $\hat{\Omega}^{(j)}$. We also assume that the incident wave $u^\mathrm{in}$ allows the cylindrical expansion written as
\begin{align}
 u^\mathrm{in}(\bm{x}) = \sum_{n=-\infty}^\infty \alpha^{(i)}_n I_n(\bm{x}-\bm{x}_0^{(i)}) \quad x\in D^{(i)},
\end{align}
for $i=1,\ldots,N$ and $\alpha^{(i)}=(\alpha^{(i)}_n)_{n\in\mathbb{Z}}$. Then, the scattering matrix $S^{(i)}$ relates incoming and outgoing waves around each $\hat{\Omega}^{(i)}$ by
\begin{align}
 S^{(i)} \left(\alpha^{(i)} + \sum_{j\neq i}T^{(ij)}B^{(j)} \right) = B^{(i)}. \label{eq:multi}
\end{align}
This linear system solves the unknown vectors $B^{(i)}$ when the incident coefficients $\alpha^{(i)}$ are provided. 

The scattering matrices and relevant vectors are of infinite size; thus we have to truncate them in practical computations. We truncate the infinite series $\sum_{n=-\infty}^\infty$ into $\sum_{n=-n_\mathrm{tr}}^{n_\mathrm{tr}}$ using Rokhlin's empirical formula \cite{coifman1993fast} given by
\begin{align}
 n_\mathrm{tr} = kd+8\log(kd+\pi),
\end{align}
where $d$ denotes the minimum distance between two centers, i.e., $d=\min_{i\neq j} |\bm{x}_0^{(i)}-\bm{x}_0^{(j)}|$.

\subsection{Scattering through periodic obstacles}
\subsubsection{Scattering matrix method}
\begin{figure}[h]
  \begin{center}
  \includegraphics[scale=0.4]{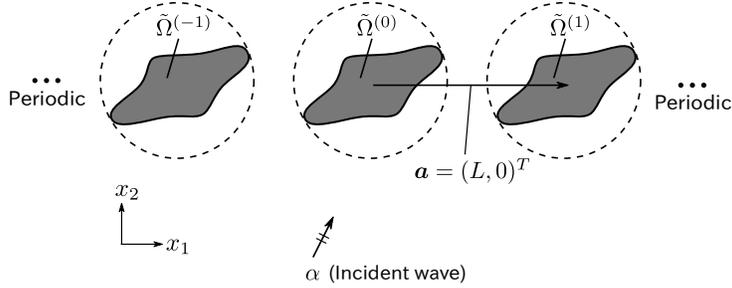}
   \caption{Scattering through a grating comprising periodic scatterers $\hat{\Omega}^{(i)}$ ($i\in\mathbb{Z}$) placed in $\mathbb{R}^2$}
   \label{fig:periodic}
  \end{center}
\end{figure}
Next, we describe a scattering matrix method for wave scattering through periodic obstacles, which was originally proposed in \cite{nicorovici1995photonic} for circular rods. As shown in \cref{fig:periodic}, we consider scatterers periodically embedded along a line in $\mathbb{R}^2$. In this study, the lattice vector $\bm{a}\in\mathbb{R}^2$ is given by $\bm{a}=(L,0)^T$ without the loss of generality, where $L>0$ denotes a given constant. We also assume that all scatterers are identical so that $\hat{\Omega}^{(i)}$ ($i\in\mathbb{Z}$) has the same scattering matrix $S$.

The scattering matrix reduces the periodic scattering problem into a system of linear algebraic equations involving the outgoing multipole coefficients $B^{(i)}$ and incident coefficients $\alpha^{(i)}$ associated with each scatterer $\hat{\Omega}^{(i)}$. To investigate this, we formally use \cref{eq:multi} to obtain
\begin{align}
 S \left(\alpha^{(i)} + \sum_{j\neq i,j=-\infty}^\infty T^{(ij)}B^{(j)} \right) = B^{(i)}, \label{eq:tmp1}
\end{align}
for each $i\in\mathbb{Z}$. We assume that the incident wave $u^\mathrm{in}$ has the quasiperiodicity $u^\mathrm{in}(\bm{x}+\bm{a})=u^\mathrm{in}(\bm{x})\e^{\ione\beta}$ for a Floquet wavenumber $\beta\in\mathbb{C}$, which is equivalent to $\alpha^{(i+1)}=\alpha^{(i)}\e^{\ione\beta}$. Then, from the Bloch--Floquet theorem, $B^{(i)}$ should satisfy the same quasiperiodic condition $B^{(i+1)}=B^{(i)}\e^{\ione\beta}$. Substituting these conditions into \cref{eq:tmp1}, we obtain
\begin{align}
 (I-S(\omega)T^\mathrm{G}(\omega,\beta))B(\omega,\beta) = S(\omega)\alpha(\omega,\beta), \label{eq:periodic}
\end{align}
where $B=B^{(0)}$ and $\alpha=\alpha^{(0)}$. The matrix $T^\mathrm{G}$ is defined by the lattice sum as follows:
\begin{align}
 T^\mathrm{G}_{ij} = \sum_{n\in\mathbb{Z}\setminus\{0\}} T^{(0n)}_{ij}\e^{\ione n\beta} = \sum_{n\in\mathbb{Z}\setminus\{0\}} O_{j-i}(-n\bm{a})\e^{\ione n\beta}. \label{eq:schlomilch}
\end{align}
This lattice sum is called a \textit{Schl\"{o}milch series} and slowly convergent if $\beta\in\mathbb{R}$ and $\omega>0$ \cite{linton2006schlomilch}. 

{Note that the proposed method is closely related to BEM with quasi-periodic Green's function \cite{porter1999rayleigh,otani2008fmm,isakari2012calderon,nose2014calculation,misawa2016fmm}. Although these approaches are more straightforward, the scattering matrix formulation is more convenient for evaluating the topological derivative, introduced in \cref{ss:td}.}

\subsubsection{Integral representation of the Schl\"{o}milch series}
Although the Schl\"{o}milch series \cref{eq:schlomilch} is convergent {for real $\omega$ and $\beta$}, we need a more rapidly convergent representation to evaluate it numerically. Moreover, we wish to establish a representation that is valid even for complex $\beta$ and $\omega$ to compute resonant frequencies $\omega$ and wavenumbers $\beta$ because they lie in the complex planes. 

\begin{figure}[h]
  \begin{center}
  \includegraphics[scale=0.3]{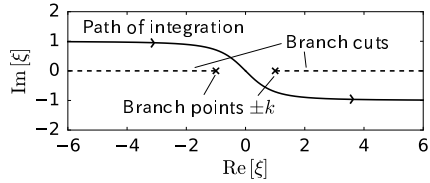}
   \caption{Path of integration for the rightmost term in the RHS of \cref{eq:OF} and branch cuts of $R(\xi)$ for $k=1$}
   \label{fig:path}
\end{center}
\end{figure}
{First, we assume that $\omega>0$ and $\beta\in\mathbb{R}$.} According to \cite{otani2008fmm}, $T^\mathrm{G}$ has the following integral representation:
\begin{align}
 T^\mathrm{G}_{ij} &= \sum_{n=1}^{s-1} O_{j-i}(-n\bm{a})\e^{\ione n\beta}
 + \sum_{n=-s+1}^{-1} O_{j-i}(-n\bm{a})\e^{\ione n\beta} + \frac{1}{\pi\ione k^{i-j}}\int_{-\infty}^\infty f_{i-j}(\xi)\dd\xi, \label{eq:OF}
 \\
 f_n(\xi) &= \frac{\e^{s(-\ione\beta - R(\xi)L)} (\xi-R(\xi))^n}{R(\xi)(1-\e^{-\ione\beta - R(\xi)L})} + \frac{\e^{s(\ione\beta - R(\xi)L)} (\xi+R(\xi))^n}{R(\xi)(1-\e^{\ione\beta - R(\xi)L})},
 \\
 R(\xi) &= {\sqrt{\xi^2-k^2}} , 
\end{align}
where the integer $s\geq 2$ is arbitrary. {This integral representation is a modified version of Linton's integral form \cite{linton2006schlomilch}.} To obtain the convergence of the integral in \cref{eq:OF}, we have to determine the branch cuts of the integrand, choose an appropriate sheet, and deform the integration path to circumvent the branch cuts. A possible choice is $R(\xi)={\sqrt{|\xi^2-k^2|}\exp(\ione\mathrm{Arg}\,(\xi^2-k^2)/2)}$, where $\mathrm{Arg}:\mathbb{C}\to(-\pi,\pi]$ is the principal argument. Here, we can choose a path of integration as the steepest descent path of $\exp(-s R(\xi)L)$ to obtain a rapid convergence. Thus, we use the path given by $\xi(t) = \pm Q(t)$ for $t\in[0,\infty)$, where $Q(t)=\sqrt{|t^2-2\ione kt|}\exp(\ione\mathrm{Arg}\,(t^2-2\ione kt)/2)$. \cref{fig:path} illustrates the path of integration and branch cuts from which we confirm that the path does not cross the branch cuts. Finally, we obtain
\begin{align}
 \int_{-\infty}^\infty f_{i-j}(\xi)\dd\xi = \int_0^\infty 
 \biggl[
 f_{i-j}(Q(t)) 
 + f_{i-j}(-Q(t))
 \biggr]
 \frac{t-\ione k}{Q(t)} \dd t.
 \label{eq:OF_2}
\end{align}
Because the integrand in \cref{eq:OF_2} is oscillatory and has a weak singularity of order $t^{-1/2}$ at $t=0$, we further apply the double-exponential formula \cite{ooura1999robust} to this integral in the practical computation.

{Although the integral expression \cref{eq:OF} is originally proposed for real $\omega$ and $\beta$, the convergence of the Fourier integral \cref{eq:OF_2} is still achieved for complex $\beta$ and $\omega$. To see this, let us evaluate the integrand in \cref{eq:OF_2} as follows:
\begin{align}
    &\biggl[
      f_{n}(Q(t)) 
      + f_{n}(-Q(t))
      \biggr]
      \frac{t-\ione k}{Q(t)}
      \notag\\
      =& \frac{1}{\sqrt{t^2-2\ione kt}} \biggl[
      \frac{(\sqrt{t^2-2\ione kt}-t+\ione kt))^n \e^{-Lst+\ione s(kL-\beta)}}{1-\e^{-Lt+\ione (kL-\beta)}}
      \notag\\
      &\hspace{50pt}+
      \frac{(\sqrt{t^2-2\ione kt}+t-\ione kt))^n \e^{-Lst+\ione s(kL+\beta)}}{1-\e^{-Lt+\ione (kL+\beta)}}
      \biggr]
      \notag\\
      =& O(t^{|n|-1}\e^{-sLt}), \quad t\to+\infty.
\end{align}
This estimation shows that the integral \cref{eq:OF_2} is convergent even if $\Im[\beta] \neq 0$ or $\Im[\omega] \neq 0$.
}
 The representation {\cref{eq:OF_2}} has an infinite number of branch points at $\beta=2m\pi\pm kL$ for $m\in\mathbb{Z}$, yielding the Rayleigh anomaly. More careful investigations \cite{nose2014calculation} show that all branch cuts of $\beta\mapsto T^\mathrm{G}$ for a fixed $\omega$ are written as $\{kL+2m\pi+\ione y \mid y\geq0,\, m\in\mathbb{Z}\}$ and $\{-kL+2m\pi-\ione y \mid y\geq0,\, m\in\mathbb{Z}\}$.

\begin{figure}[h]
  \begin{center}
  \includegraphics[scale=0.7]{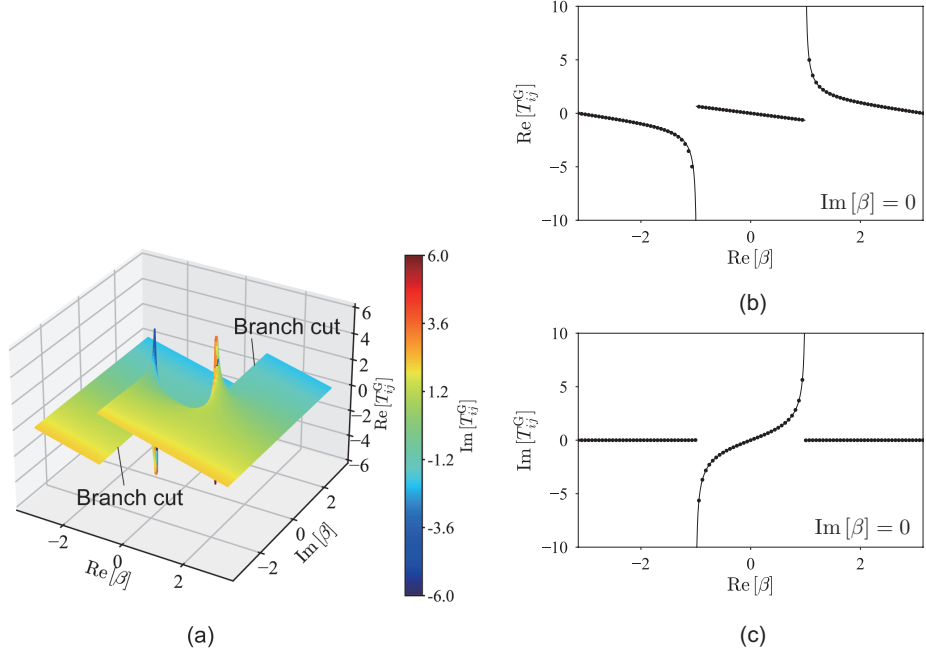}
   \caption{{Schl\"{o}milch series $T^\mathrm{G}_{ij}$ for $i-j=1$, $k=1$, and $L=1$. (a) Values of $T^\mathrm{G}_{ij}$ computed using the integral representation \cref{eq:OF_2}. (b) and (c) Real and imaginary parts of $T^\mathrm{G}_{ij}$ for real $\beta$, respectively. The values are calculated using the lattice sum \cref{eq:schlomilch} (dots) and integral representation \cref{eq:OF_2} (solid lines).}}
   \label{fig:schlomilch}
\end{center}
\end{figure}
{In \cref{fig:schlomilch}, we plot the values of $T^\mathrm{G}_{ij}$ calculated using the lattice sum \cref{eq:schlomilch} and integral representation \cref{eq:OF_2} with $s=2$ to validate the expressions. The lattice sum \cref{eq:schlomilch} is truncated at $|n| = 10^8$. In this computation, the comparison between \cref{eq:schlomilch} and \cref{eq:OF_2} is given only along the real axis because the lattice sum \cref{eq:schlomilch} is divergent otherwise. The results show that $T^\mathrm{G}_{ij}$ is smoothly extended into the complex $\beta$-plane except for the branch cuts. In addition, the values are in good agreement with the truncated lattice sum. 
}

\subsection{Modal analysis}\label{ss:guided}
The scattering-matrix formalism \cref{eq:periodic} allows us to perform guided- and leaky-mode analysis by finding pairs $(\omega,\beta)$ such that the linear system \cref{eq:periodic} has a nontrivial solution $B$ without any incident field $\alpha$. This is a nonlinear eigenvalue problem for the matrix-valued function $I-ST^\mathrm{G}$ when either $\omega$ or $\beta$ is fixed in $\mathbb{C}$. Therefore, it can be solved using a gradient- or contour integral-based algorithm. In this study, we adopt the Sakurai--Sugiura method (SSM) \cite{asakura2009numerical}, which determines an eigenpair $(z,\phi)$ of $F(z)\phi=0$, where $F$ denotes a matrix-valued and possibly nonlinear function, within a closed path $C$ in $\mathbb{C}$ by integrating $u^H F^{-1}v$ for some $u$ and $v$ on $C$ and converting the nonlinear eigenvalue problem into a generalized eigenvalue problem. {The SSM can find multiple eigenvalues (even if they are degenerated) in $C$ when an appropriate parameter is given in the algorithm.} {This approach is originally proposed and validated by Nose and Nishimura \cite{nose2014calculation} with a fast multipole method. They applied the SSM to find nonlinear eigenvalues of a coefficient matrix that arises in a BEM with quasi-periodic Green's function.} We refer to \cite{asakura2009numerical} for more details about the SSM algorithm.

 {For a fixed $\omega>0$, a guided mode propagates along periodic obstacles without attenuation in space, meaning that a resonant wavenumber $\beta$ is real. On the other hand, if $\beta$ is complex, the corresponding mode decays exponentially as it travels along the structure. In this case, we say that a resonant mode is leaky.} A leaky mode $u$ satisfies the original boundary value problem \cref{eq:bvp_1,eq:bvp_2,eq:bvp_3,eq:bvp_4,eq:bvp_5}. Furthermore, if the pair {$\omega>0$ and $\beta\in\mathbb{R}$} lies in the radiation continuum, i.e., $\omega^2/c^2-(\beta+2n\pi)^2/L^2 > 0$ for some $n\in\mathbb{Z}$, then the bound state is called a BIC.

\section{Topology optimization}
In this section, we design the shape and topology of a unit structure comprising a periodic waveguide such that it exhibits desirable resonant properties. To this end, we use a topology optimization approach \cite{bendsoe2013topology} to seek an optimal material distribution for a given objective functional. Here, the objective functional is set as {$(\Im[\beta])^2$ with} a resonant wavenumber {$\beta$}. {To apply topology optimization, we need a sensitivity of the given objective functional with respect to a geometrical perturbation, called a topological derivative. In this section, we first derive a novel expression of the topological derivative for the eigenvalue problem in \cref{ss:td} and then explain the algorithm for the topology optimization in \cref{ss:levelset}.}

\subsection{Topological derivative}\label{ss:td}
For an effective optimization algorithm, we need sensitivity with respect to a small perturbation of the geometry of a unit structure. In this subsection, we derive topological derivatives \cite{sokolowski1999topological} related to resonant properties of the periodic waveguide. 

\subsubsection{Scattering matrix}
\begin{figure}[h]
  \begin{center}
  \includegraphics[scale=0.32]{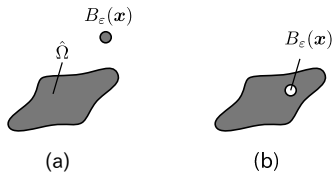}
   \caption{Topological change around a single scatterer $\hat{\Omega}$. (a) Case that a small disk $B_\eps$ appears in the exterior $\Omega$. (b) Case that a small disk $B_\eps$ appears in the interior $\hat{\Omega}$.}
   \label{fig:td}
\end{center}
\end{figure}
Here, we first investigate the perturbation of the scattering matrix $S$ associated with a single scatterer $\hat{\Omega}$ by a small particle added at a point $x$ in either $\Omega=\mathbb{R}^2\setminus\overline{\hat{\Omega}}$ or $\hat{\Omega}$. Let $B_\eps(\bm{x})$ be an open disk of radius $\eps$ centered at $x$. First, we consider the case of $x\in \Omega$ and assume that the disk $B_\eps(\bm{x})$ is characterized by $\hat{\rho}$ and $\hat{\kappa}$, i.e. $B_\eps(\bm{x})$ comprises the same material filling in $\hat{\Omega}$, as shown in \cref{fig:td} (a). For sufficiently small radius $\eps>0$, let $\delta S$ denote the perturbation of $S$, i.e. 
\begin{align}
 S(\omega;\hat{\Omega}\cup B_\eps(\bm{x})) = S(\omega;\hat{\Omega}) + \delta S. 
\end{align}
Recall that the scattering matrix $S$ is given by
\begin{align}
 S_{nn^\prime} = \frac{\ione(-1)^n}{4}\langle I_{-n}, u_{n^\prime}  \rangle_{\partial\hat{\Omega}},
\end{align}
which yields the variation
\begin{align}
 \delta S_{nn^\prime} = \frac{\ione(-1)^n}{4}\langle I_{-n}, \delta u_{n^\prime}  \rangle_{\partial\hat{\Omega}} + \frac{\ione(-1)^n}{4}\langle I_{-n}, u_{n^\prime} + \delta u_{n^\prime}  \rangle_{\partial B_\eps(\bm{x})}, \label{eq:delta_S}
\end{align}
where $u_{n^\prime}$ denotes the solution of the boundary value problem \cref{eq:bvp_1,eq:bvp_2,eq:bvp_3,eq:bvp_4,eq:bvp_5} for $u^\mathrm{in}(\bm{x})=I_{n^\prime}(\bm{x}-\bm{x}_0)$, and $u_{n^\prime}+\delta u_{n^\prime}$ represents the solution of the boundary value problem defined by replacing $\Omega$ with $\Omega\setminus\overline{B_\eps(\bm{x})}$ in \cref{eq:bvp_1,eq:bvp_2,eq:bvp_3,eq:bvp_4,eq:bvp_5}.

Let $\tilde{u}_n$ be an adjoint variable satisfying the Helmholtz equations
\begin{align}
 \nabla^2 \tilde{u}_n(\bm{x}) + \frac{\omega^2}{c^2} \tilde{u}_n(\bm{x}) = 0 &\quad x\in \Omega,
 \\
 \nabla^2 \tilde{u}_n(\bm{x}) + \frac{\omega^2}{\hat{c}^2} \tilde{u}_n(\bm{x}) = 0 &\quad x\in\hat{\Omega},
\end{align}
 and the Sommerfeld radiation condition. Then, the reciprocity theorem yields
\begin{align}
 \int_{\partial\hat{\Omega}\cup\partial B_\eps(\bm{x})} \left(
 \tilde{u}_n|_+ \frac{\partial \delta u_{n^\prime}} {\partial n}\Big|_+ - \delta u_{n^\prime}|_+ \frac{\partial \tilde{u}_n}{\partial n}\Big|_+
 \right)
 \dd\Gamma = 0, \label{eq:recipro_1}
 \\ 
 \int_{\partial\hat{\Omega}} \left(
 \tilde{u}_n|_- \frac{\partial \delta u_{n^\prime}} {\partial n}\Big|_- - \delta u_{n^\prime}|_- \frac{\partial \tilde{u}_n}{\partial n}\Big|_-
 \right)
 \dd\Gamma = 0. \label{eq:recipro_2}
\end{align}
From \cref{eq:recipro_1,eq:recipro_2}, we have
\begin{align}
 &\quad \langle I_{-n}, \delta u_{n^\prime}  \rangle_{\partial\hat{\Omega}} 
 \notag\\
 &= 
 \langle I_{-n}, \delta u_{n^\prime}  \rangle_{\partial\hat{\Omega}}
 - \frac{1}{\rho} \int_{\partial\hat{\Omega}\cup\partial B_\eps(\bm{x})} \left(
 \tilde{u}_n|_+ \frac{\partial \delta u_{n^\prime}} {\partial n}\Big|_+ - \delta u_{n^\prime}|_+ \frac{\partial \tilde{u}_n}{\partial n}\Big|_+
 \right)
 \dd\Gamma
 \notag\\
 &= 
 \int_{\partial\hat{\Omega}} \left[\frac{1}{\rho}\frac{\partial \delta u_{n^\prime}}{\partial n}\Big|_+ \left(-\tilde{u}_n|_+ + \tilde{u}_n|_- -\rho I_{-n} \right) - \delta u_{n^\prime}|_+ \left(-\frac{1}{\rho}\frac{\partial \tilde{u}_n}{\partial n}\Big|_+ + \frac{1}{\hat{\rho}}\frac{\partial \tilde{u}_n}{\partial n}\Big|_- - \frac{\partial I_{-n}}{\partial n} \right)  \right]
 \notag\\
 &\quad + \frac{1}{\rho} \langle \tilde{u}_n, \delta u_{n^\prime}  \rangle_{\partial B_\eps(\bm{x})}.
 \label{eq:tmp10}
\end{align}
Imposing the boundary conditions
\begin{align}
 \tilde{u}_n|_+ = \tilde{u}_n|_- - \rho I_{-n} \quad \bm{x}\in\partial\hat{\Omega}, \label{eq:bc_adjoint_1}
 \\
 \frac{1}{\rho}\frac{\partial \tilde{u}_n}{\partial n}\Big|_+ = \frac{1}{\hat{\rho}}\frac{\partial \tilde{u}_n}{\partial n}\Big|_- - \frac{\partial I_{-n}}{\partial n} \quad \bm{x}\in\partial\hat{\Omega}, \label{eq:bc_adjoint_2}
\end{align}
\cref{eq:tmp10} is reduced to $\langle I_{-n}, \delta u_{n^\prime}  \rangle_{\partial\hat{\Omega}} = \frac{1}{\rho} \langle \tilde{u}_n, \delta u_{n^\prime}  \rangle_{\partial B_\eps(\bm{x})}$. 
Substituting this into \cref{eq:delta_S}, we obtain
\begin{align}
 \delta S_{nn^\prime} &= \frac{\ione(-1)^n}{4} \left( \frac{1}{\rho}\langle \tilde{u}_{n},\delta u_{n^\prime} \rangle_{\partial B_\eps(\bm{x})}  + \langle I_{-n}, u_{n^\prime} + \delta u_{n^\prime}  \rangle_{\partial B_\eps(\bm{x})} \right),
 \notag\\
 &= \frac{\ione(-1)^n}{4} \left(  
 \langle I_{-n}, u_{n^\prime} \rangle_{\partial B_\eps(\bm{x})} + \langle \frac{1}{\rho}\tilde{u}_{n}+I_{-n},\delta u_{n^\prime} \rangle_{\partial B_\eps(\bm{x})}
 \right),
 \notag\\
 &= \frac{\ione(-1)^n}{4} \langle \frac{1}{\rho}\tilde{u}_{n}+I_{-n},\delta u_{n^\prime} \rangle_{\partial B_\eps(\bm{x})}. \label{eq:tmp_111}
\end{align}
 Here, we have used the reciprocity $\langle I_{-n}, u_{n^\prime} \rangle_{\partial B_\eps(\bm{x})}=0$. From the boundary conditions \cref{eq:bc_adjoint_1,eq:bc_adjoint_2}, we have $\frac{1}{\rho}\tilde{u}_n + I_{-n} = u_{-n}$. This formula further simplifies \cref{eq:tmp_111} as follows:
 \begin{align}
  \delta S_{nn^\prime} &= \frac{\ione(-1)^n}{4} \langle u_{-n},\delta u_{n^\prime} \rangle_{\partial B_\eps(\bm{x})}. \label{eq:delta_S_2}
 \end{align}

 We can no longer simplify the expression \cref{eq:delta_S_2}. However, we are only interested in the asymptotic behavior of $\delta S$ for $\eps\to 0$. This can be achieved by expanding $\delta u_{n^\prime}$ with respect to $\eps$ around the point $x$.
 According to \cite{nakamoto2017levelsetbased}, we have
\begin{align}
 & \langle u_{-n},\delta u_{n^\prime} \rangle_{\partial B_\eps(\bm{x})} 
 \notag\\
 =& 
 \pi \eps^2
 \Biggl[
 \frac{2(\hat{\rho}-\rho)}{\rho+\hat{\rho}} \nabla u_{-n}(\bm{x})\cdot\nabla u_{n^\prime}(\bm{x}) 
 + \omega^2 \rho\left(\frac{1}{\hat{\kappa}} - \frac{1}{\kappa} \right) u_{-n}(\bm{x}) u_{n^\prime}(\bm{x})
 \Biggr] + O(\eps^3). \label{eq:asymp}
\end{align}

Now, we define a topological derivative of $f$, denoted by $\mathcal{D}_\mathrm{T} f$, as
\begin{align}
 \mathcal{D}_\mathrm{T} f = \lim_{\eps\to0} \frac{\delta f}{\pi \eps^2},
\end{align}
where $\delta f$ denotes the variation of $f$ due to the topological change. Then, \cref{eq:asymp} gives the final expression for the topological derivative as follows:
\begin{align}
 \mathcal{D}_\mathrm{T} S_{nn^\prime} = \frac{\ione(-1)^n}{4}
 \Biggl[
 \frac{2(\hat{\rho}-\rho)}{\rho+\hat{\rho}} \nabla u_{-n}(\bm{x})\cdot\nabla u_{n^\prime}(\bm{x}) 
 + \omega^2 \rho \left( \frac{1}{\hat{\kappa}} - \frac{1}{\kappa} \right) u_{-n}(\bm{x}) u_{n^\prime}(\bm{x})
 \Biggr]
 \notag\\
 \bm{x}\in \Omega.
\end{align}

We can treat the case of $x\in\hat{\Omega}$ (\cref{fig:td} (b)) in the same manner. In this case, the topological derivative is given by
\begin{align}
 \mathcal{D}_\mathrm{T} S_{nn^\prime} = \frac{\ione(-1)^n}{4}
 \Biggl[
 \frac{2(\rho-\hat{\rho})}{\rho+\hat{\rho}} \nabla u_{-n}(\bm{x})\cdot\nabla u_{n^\prime}(\bm{x}) 
 + \omega^2 \rho \left(\frac{1}{\kappa}- \frac{1}{\hat{\kappa}} \right) u_{-n}(\bm{x}) u_{n^\prime}(\bm{x})
 \Biggr]
 \notag\\
 \bm{x}\in \hat{\Omega}.
\end{align}

\subsubsection{Resonant wavenumber}
The topological perturbation changes the distribution of the resonant frequencies and wavenumbers of the periodic system, characterized by the equation \cref{eq:periodic} with $\alpha=0$. We fix $\omega$ and investigate the variation in $\beta$ caused by the topological change.

Suppose that the equation
\begin{align}
 (I-S(\omega;\Omega)T^\mathrm{G}(\omega,\beta(\Omega)))B(\omega,\beta(\Omega)) = 0 \label{eq:tmp30}
\end{align}
and its perturbed system
\begin{align}
 (I-S(\omega;\hat{\Omega}\cup B_\eps(\bm{x}))T^\mathrm{G}(\omega,\beta(\hat{\Omega}\cup B_\eps(\bm{x}))))B(\omega,\beta(\hat{\Omega}\cup B_\eps(\bm{x}))) = 0 \label{eq:tmp31}
\end{align}
have nontrivial solutions. We evaluate the difference $\delta \beta = \beta(\hat{\Omega}\cup B_\eps(\bm{x})) - \beta(\Omega)$. From \cref{eq:tmp30,eq:tmp31}, we have
\begin{align}
 S\frac{\partial T^\mathrm{G}}{\partial \beta} B \delta \beta = -\delta S T^\mathrm{G}B + (I-ST^\mathrm{G})\delta B - \delta S(\delta T^\mathrm{G}B + (T^\mathrm{G}+\delta T^\mathrm{G})\delta B) - S\delta T^\mathrm{G} \delta B. \label{eq:tmp32}
\end{align}
Let $\tilde{B}$ be a left eigenvector satisfying $\tilde{B}^H(I-ST^\mathrm{G})=0$. Then, multiplying both sides of \cref{eq:tmp32} by $\tilde{B}^H$, we have
\begin{align}
 \tilde{B}^H S\frac{\partial T^\mathrm{G}}{\partial \beta} B \delta \beta = -\tilde{B}^H \delta S T^\mathrm{G}B - \tilde{B}^H \delta S(\delta T^\mathrm{G}B + (T^\mathrm{G}+\delta T^\mathrm{G})\delta B) - \tilde{B}^HS\delta T^\mathrm{G} \delta B,
\end{align}
which gives the topological derivative
\begin{align}
 \mathcal{D}_\mathrm{T}\beta = -\frac{\tilde{B}^H (\mathcal{D}_\mathrm{T} S) T^\mathrm{G} B}{\tilde{B}^H S \frac{\partial T^\mathrm{G}}{\partial \beta} B}.
\end{align}
Here, we have neglected the higher-order variations.

\subsection{Algorithm for the topology optimization}\label{ss:levelset}
\begin{figure}[h]
  \begin{center}
  \includegraphics[scale=0.4]{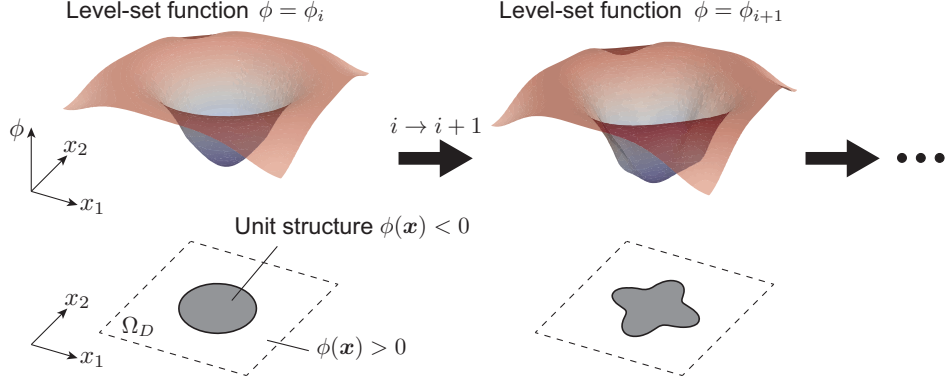}
   \caption{{Schematic illustration of the level-set-based topology optimization method.}}
   \label{fig:levelset}
\end{center}
\end{figure}
Herein, we perform the topology optimization to {find a shape of a unit structure $\hat{\Omega}$ that minimizes the objective functional $J=(\Im[\beta])^2$ for a fixed $\omega$. If the objective value attains $J=0$, the obtained shape $\hat{\Omega}$ should exhibit a BIC at the target frequency. To this end, we employ a level-set-based topology optimization algorithm.}. First, we define a scalar function $\phi:\Omega_D\to \mathbb{R}$, called a \textit{level-set function}, within a fixed design domain $\Omega_D\subset \mathbb{R}^2$. The level-set function $\phi$ gives the material distribution in $\Omega_D$ by
\begin{align}
 \hat{\Omega} = \{ x\in \Omega_D \mid \phi(\bm{x}) < 0  \},
 \\
 \Omega_D \setminus \overline{\hat{\Omega}} = \{ x\in \Omega_D \mid \phi(\bm{x}) > 0  \},
 \\
 \partial\hat{\Omega} = \{ x\in \Omega_D \mid \phi(\bm{x}) = 0  \}.
\end{align}
{Instead of seeking an optimal shape of $\hat{\Omega}$ directly, level-set-based topology optimization methods optimize the distribution of $\phi$ using iterative algorithms. This procedure is illustrated in \cref{fig:levelset}.} Following \cite{amstutz2006new}, we iteratively update the level-set function $\phi$ by the following formula:
\begin{align}
 \phi_{i+1}(\bm{x}) = (1-\Delta_i({\mathcal{T}_i},{\phi_i})_{L^2(\Omega_D)}) \phi_i(\bm{x}) + \Delta_i {\mathcal{T}_i(\bm{x})}, \label{eq:amstutz}
\end{align}
where $\phi_i$ denotes the level-set function at $i$th step, $\Delta_i>0$ denotes a step length, {$\mathcal{T}_i$} represents the topological derivative of {the} objective functional {$J$ corresponding to $\hat{\Omega}$ at $i$th step}, and $(\cdot,\cdot)_{L^2(\Omega_D)}$ denotes the $L^2$ inner product in $\Omega_D$ defined by
\begin{align}
 (f,g)_{L^2(\Omega_D)} = \int_{\Omega_D}f(\bm{x})g(\bm{x})\dd\Omega,
\end{align}
for scalar functions $f$ and $g$ in $\Omega_D$. In the iterative algorithm, the functions $\phi_i$ and $\mathcal{T}_i$ are discretized using the B-spline basis functions \cite{isakari2017topology}. {Once the iterative procedure \cref{eq:amstutz} reaches convergence, we terminate the algorithm and obtain the optimal shape of $\hat{\Omega}$ corresponding to $\phi_i$.}

\section{Numerical examples}
In this section, we first verify the proposed method and examine the correctness of the new topological derivative. Subsequently, we present a numerical example of the topology optimization that designs a resonant waveguide exhibiting a BIC at a given frequency.

\subsection{Verification of the scattering matrix method}
\begin{figure}[h]
  \begin{center}
  \includegraphics[scale=0.3]{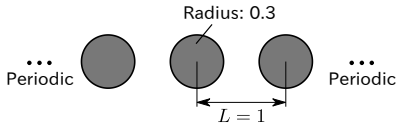}
   \caption{{Disks} placed periodically in the $x_1$ direction.}
   \label{fig:veri_problem}
\end{center}
\end{figure}
First, we verify that the proposed method determines a BIC accurately. As shown in \cref{fig:veri_problem}, we consider a waveguide comprising circular elastic materials of radius $0.3$ with mass density $\hat{\rho}=2$ and bulk modulus $\hat{\kappa}=1$ embedded in the background medium characterized by $\rho=1$ and $\kappa=1$. According to \cite{nose2014calculation}, this waveguide has a complex eigenvalue $\beta=0.591931+0.034843\ione$ for $\omega=6.2831$.

\begin{figure}[h]
  \begin{minipage}[b]{0.5\hsize}
    \centering
   \includegraphics[scale=0.3]{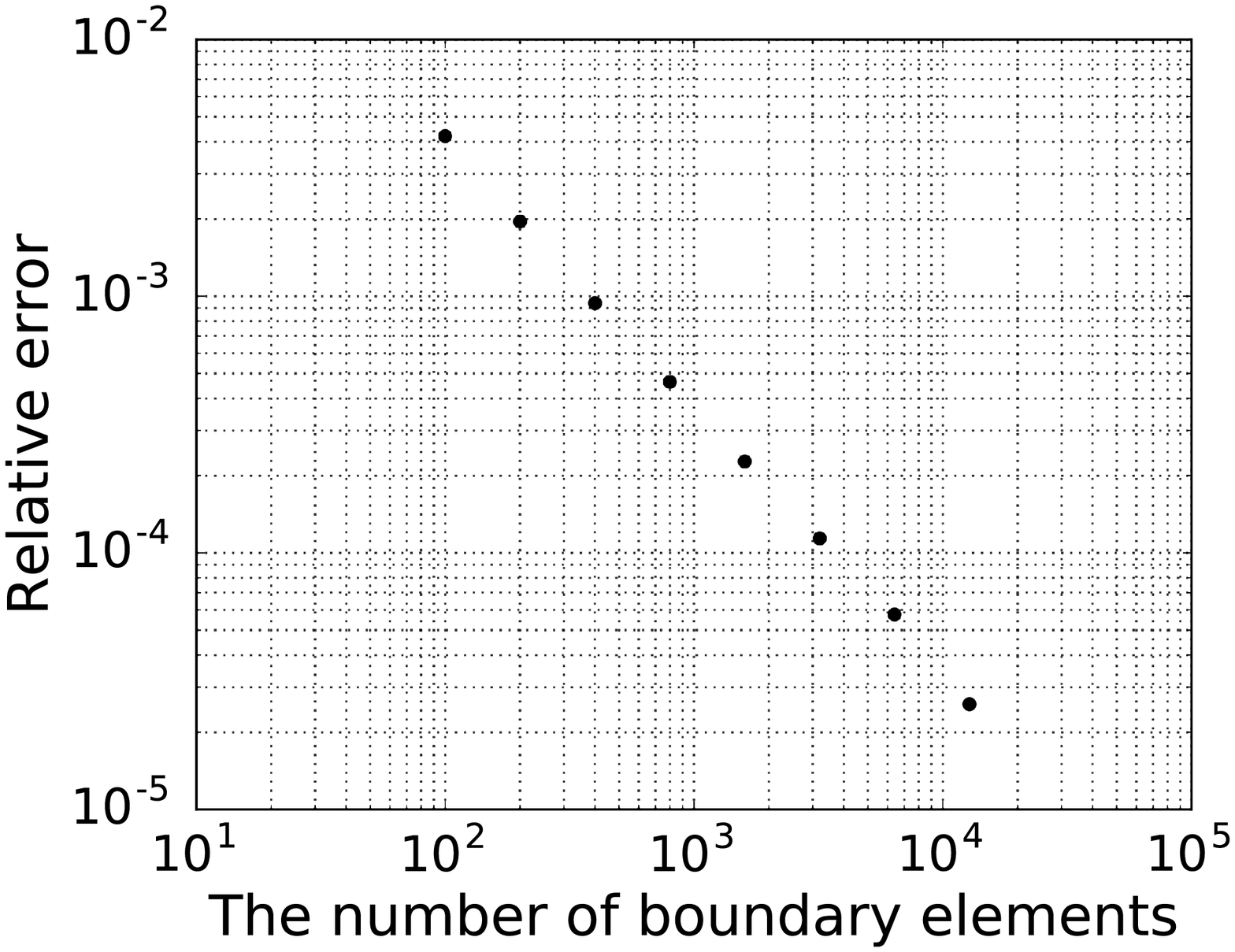}
    \subcaption{}
  \end{minipage}
  \begin{minipage}[b]{0.5\hsize}
    \centering
   \includegraphics[scale=0.3]{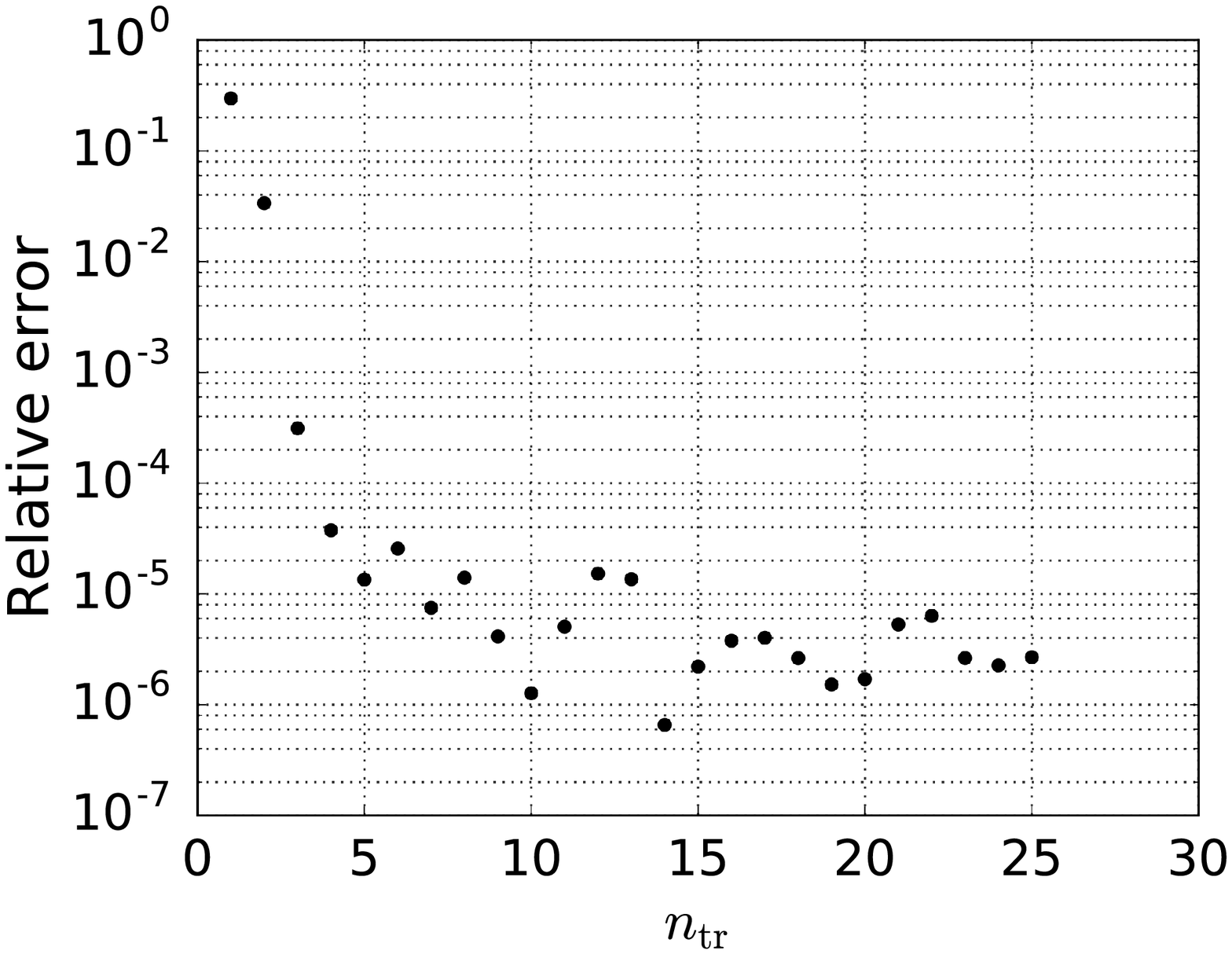}
    \subcaption{}
  \end{minipage}
 \caption{Relative error of an eigenvalue inside the path $C$. (a) The case when the number of boundary elements varies. (b) The case when the number of terms $n_\mathrm{tr}$ (defined in \Cref{ss:multi}) varies.}
 \label{fig:veri_result}
\end{figure}
We employ the SSM for a circular contour path $C$ of radius $0.4$ centered at $0.5+0.0\ione$ in the complex $\beta$-plane to determine the eigenvalue. The SSM algorithm performs contour integration along $C$ using the trapezoidal rule with 32 subintervals. Note that the path $C$ does not cross any branch cut for $\omega=6.2831$. We discretize a unit disk in \cref{fig:veri_problem} using piecewise constant boundary elements. The number of boundary elements is denoted by $N$. 

First, we perform the eigenvalue analysis for each $N=N_i$ ($i=0,1,\ldots,8$), where $N_i$ is defined as $N_i=100\times 2^i$, and fixed $n_\mathrm{tr}=20$. \cref{fig:veri_result} (a) shows the relative error of a unique eigenvalue $\beta(N)$ inside $C$ defined by $|\beta(N_{i+1}) - \beta(N_i)|/|\beta(N_i)|$ for each $N$. The result shows that the relative error decreases monotonically and converges at the rate of $O(N^{-1})$. For $N=12{,}800$, the obtained eigenvalue is $\beta=0.59259+0.035012\ione$, which is close to the value reported in \cite{nose2014calculation}. Further, we fix $N$ at $N=12{,}800$ and define a relative error in an analogous manner for $n_\mathrm{tr}=1,2,\ldots,30$. \cref{fig:veri_result} (b) shows the result of the error analysis. The result of the error analysis shows that the error monotonically decreases until it reaches around $10^{-5}$, which is close to the value at $N=12{,}800$ in \cref{fig:veri_result} (a). From these convergence tests, we conclude that the proposed method can determine resonant wavenumbers correctly.

\subsection{Topological derivative}
In this section, we examine the correctness of the topological derivative formulated in \Cref{ss:td} through a numerical experiment.

In this experiment, we use the same parameters and configuration as those used in the previous example. We compare the derivative $(\mathcal{D}_\mathrm{T}\beta)(\bm{x})$ with the finite difference $(\beta(\hat{\Omega}\cup B_\eps(\bm{x})) - \beta(\hat{\Omega}))/(\pi \eps^2)$ for some center $x$ and small radius $\eps>0$ to verify the topological derivative.

\begin{figure}[h]
  \begin{minipage}[b]{0.5\hsize}
    \centering
   \includegraphics[scale=0.28]{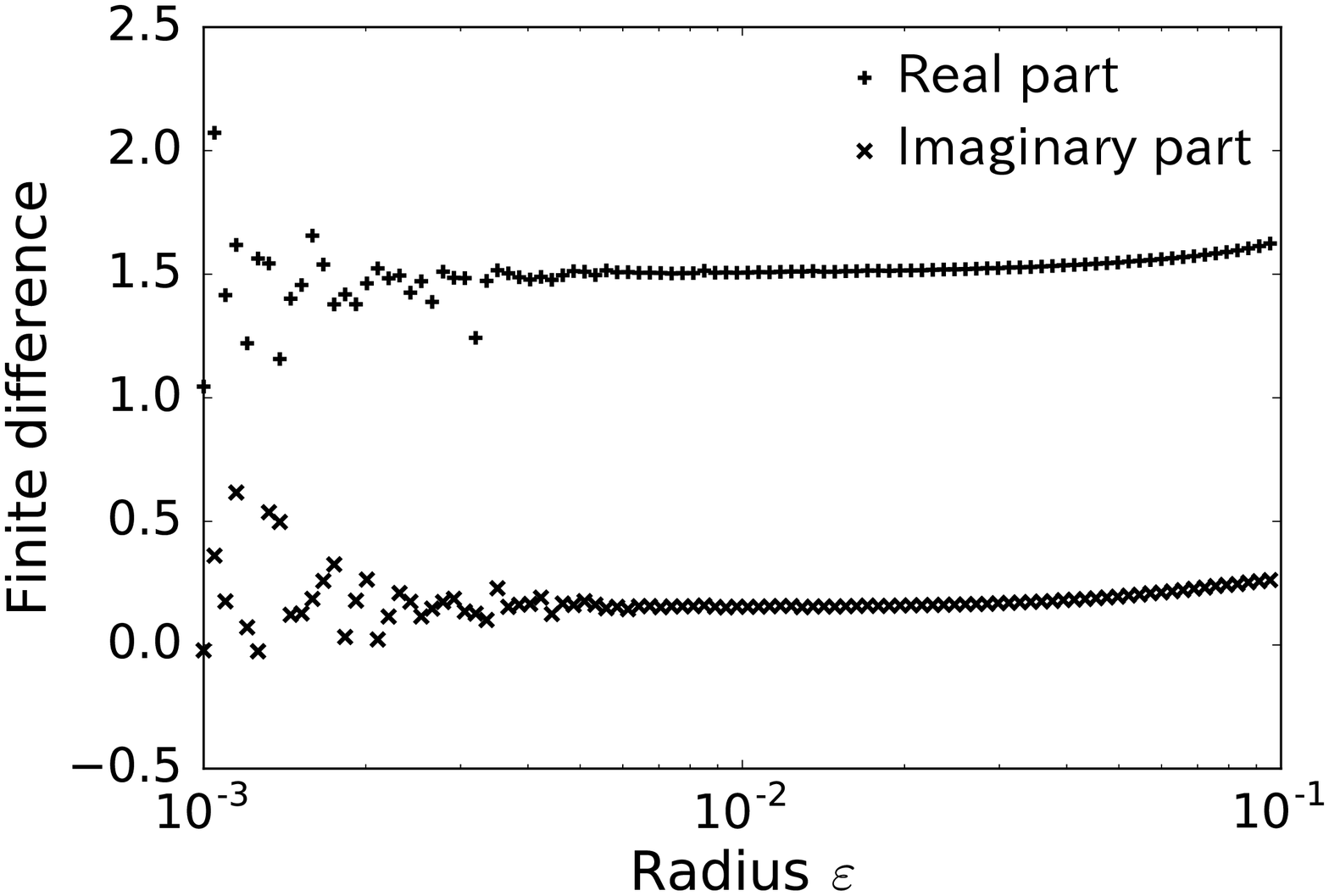}
    \subcaption{}
  \end{minipage}
  \begin{minipage}[b]{0.5\hsize}
    \centering
   \includegraphics[scale=0.28]{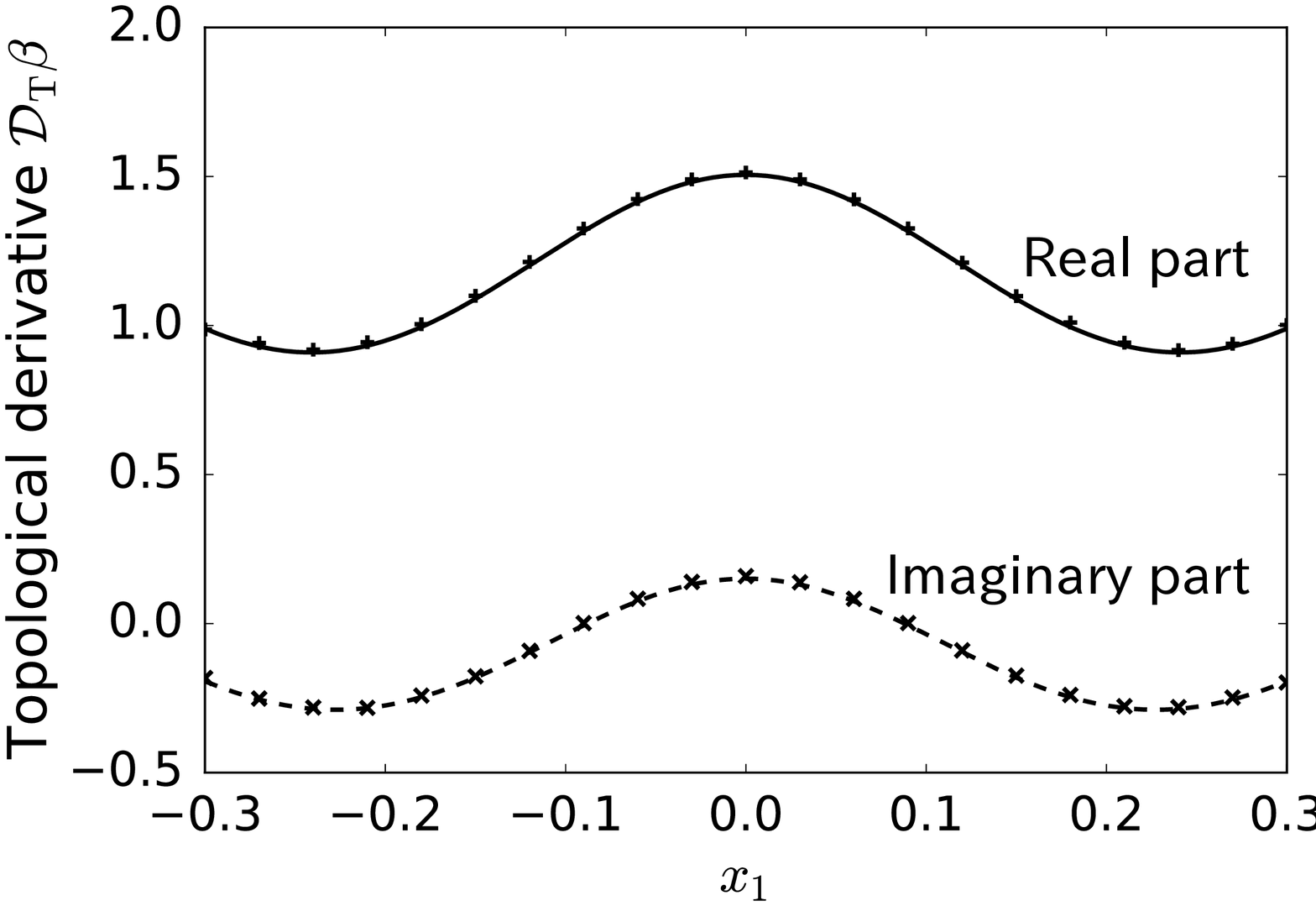}
    \subcaption{}
  \end{minipage}
 \caption{Topological derivatives $(\mathcal{D}_\mathrm{T}\beta)(\bm{x})$ and finite difference $(\beta(\hat{\Omega}\cup B_\eps(\bm{x})) - \beta(\hat{\Omega}))/(\pi \eps^2)$. (a) Finite difference versus the radius $\eps$ for $\bm{x}=(0.0,0.4)^T$. (b) Comparison between the topological derivative and finite difference for $\eps=0.02$ along the line $x_2=0.4$. The solid and dashed lines indicate the real and imaginary parts of the topological derivative, respectively. The markers express the finite difference.}
 \label{fig:td_floquet}
\end{figure}
First, we fix the center at $\bm{x}=(0.0,0.4)^T$ and investigate an appropriate radius $\eps$. \Cref{fig:td_floquet} (a) illustrates the behavior of the finite difference approximation of $\beta$ with respect to the radius $\eps$. The result shows that the approximation almost converges at $\eps=0.02$ and oscillates for smaller $\eps$ due to the loss of significant digits in computing the numerator. Thus, we can expect that $\eps=0.02$ produces a reasonable approximation to the topological derivative.

Then, we use $\eps=0.02$ and compare the approximation and topological derivative. \Cref{fig:td_floquet} (b) shows the approximation and derivative along the line $x_2=0.4$, illustrating that both values are consistent. Therefore, we conclude that the proposed topological derivative is accurate.

\subsection{Band calculation}
\begin{figure}[h]
 \centering
 \includegraphics[scale=1.0]{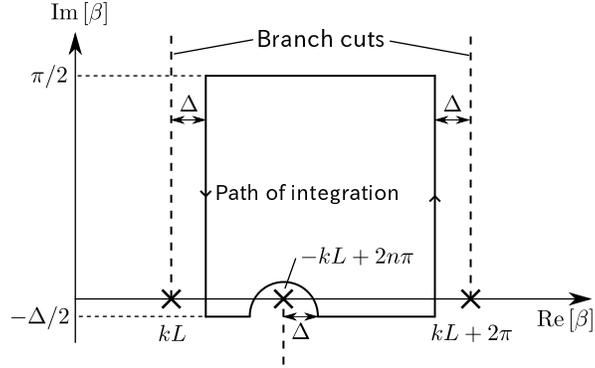}
 \caption{Path of integration used in the SSM algorithm for band calculation. The crosses denote the branch points of the function $\beta\mapsto T^\mathrm{G}(\omega,\beta)$, and the dashed lines denote the corresponding cuts.}
 \label{fig:band_path}
\end{figure}
In the previous experiments, we have focused on a single resonant wavenumber. However, we are often interested in how the eigenvalue depends on the frequency, i.e., phononic band structure. Because of the quasiperiodicity and time-reversal symmetry, it suffices to find eigenvalues in $\{\beta\in\mathbb{C} \mid kL < \mathrm{Re}\,[\beta] < kL+2\pi,\ \mathrm{Im}\,[\beta]\geq 0 \}$ for a fixed $\omega$. In the following numerical experiment, we set a path of integration $C$ for the SSM algorithm (\cref{fig:band_path}) with $\Delta=0.01$.

\begin{figure}[h]
 \centering
 \includegraphics[scale=0.25]{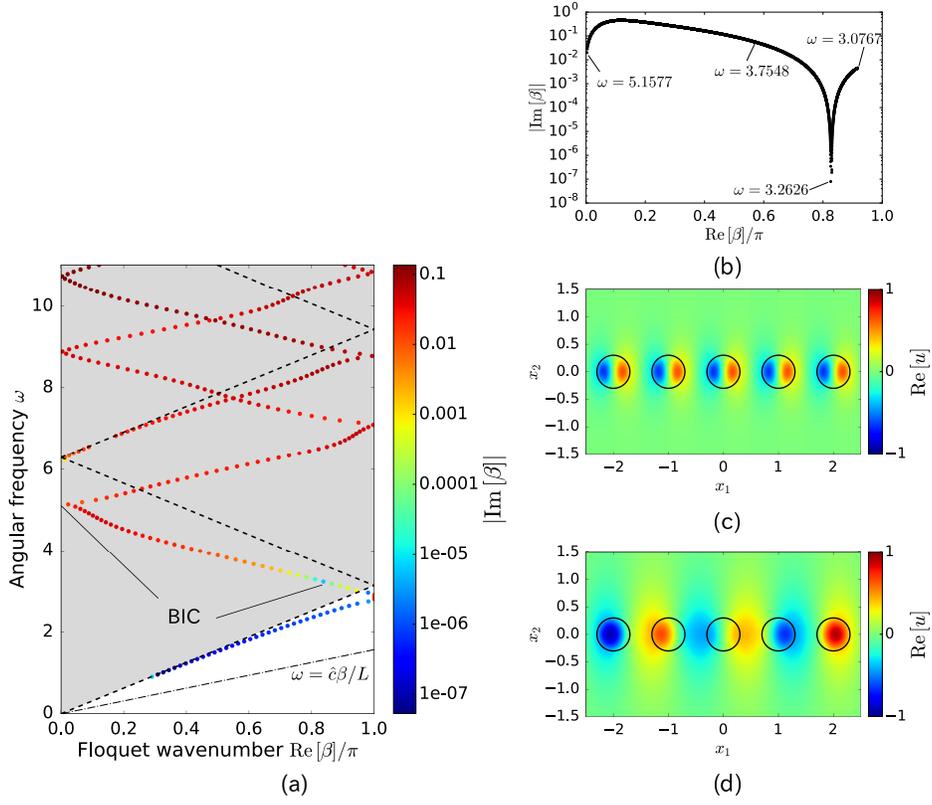}
 \caption{Result of the band calculation and BIC obtained from the analysis. (a) Band structure obtained using SSM. The dashed lines indicate the lightlines $\beta=\pm kL + 2n\pi$ with $n\in\mathbb{Z}$. (b) Imaginary part of the eigenvalues $\beta$ in $3.0767\leq \omega \leq 5.1577$. (c) Mode profile of the BIC at $(\omega,\beta)=(5.1580-1.1732\times10^{-5}\ione,0)$. (d) Mode profile of the BIC at $(\omega,\beta)=(3.2626,2.5952-7.9502\times10^{-8}\ione)$.}
 \label{fig:band_init}
\end{figure}
\Cref{fig:band_init} (a) shows the plot of the band structure obtained using the proposed method. {In the diagram, the computed eigenvalues $\beta$ are replaced with $\beta+2m\pi$, where $m$ is an integer that satisfies $0 \leq \Re[\beta+2m\pi] \leq \pi$.} {The obtained band diagram has some similar features to that of planar waveguides \cite{hu2009understanding}. For example, the diagram (\cref{fig:band_init} (a)) shows that the third band departs from the lightline at around $\omega=6.2$ (cutoff frequency). In addition, the first band $\omega(\beta)$ satisfies $\hat{c}\beta/L \leq \omega(\beta) \leq c\beta/L$. They are typical characteristics of the waveguide dispersion.} The figure shows that the eigenvalues outside the radiation continuum, which is the region below the lightlines, have small imaginary parts, thus forming guided modes along the periodic structure. Within the radiation continuum (gray-shaded region in \cref{fig:band_init} (a)), almost every eigenmode is leaky due to its nonzero imaginary part. However, we find a significantly small imaginary part within the continuum around $\omega=3.26$ and $5.16$. \cref{fig:band_init} (b) shows that the absolute values of the imaginary part decrease rapidly around $\omega=3.2626$ and $\omega=5.1577$, indicating that two BICs exist around the points. The latter point stands for a symmetry-protected BIC \cite{hsu2016bound} because it lies on the $\Gamma$ point ($\beta=0$). {As long as the parity symmetry with respect to $x_1\to -x_1$ is preserved and the material parameters satisfy a certain condition, there exists at least one symmetry-protected BIC with $\Re[\beta]=0$ \cite{bonnet-bendhia1994guided,shipman2010resonant}.} Further, we conducted an eigenvalue analysis to find a resonant $\omega$ for fixed $\beta=0$. We obtained that $(\omega,\beta)=(5.1580-1.1732\times10^{-5}\ione,0)$ is an eigenpair, whose mode profile is illustrated in \cref{fig:band_init} (c). \cref{fig:band_init} (d) shows the resonant mode corresponding to $(\omega,\beta)=(3.2626,2.5952-7.9502\times10^{-8}\ione)$. From the mode profiles, we observe that the fields are strongly confined around the structure without radiation. {This type of BICs on the second band with $\Re[\beta]\neq 0$ are already reported and discussed for circular inclusions \cite{bulgakov2014bloch,yuan2018bound}.}

\subsection{Topology optimization}
From the previous subsection, we observed that the periodic array of circular cylinders exhibits some BICs. 
{Although only the two BICs are found in the band diagram, the existence of BICs in a higher frequency regime is reported for a simple geometry \cite{bonnet-bendhia1994guided}.} In this section, we show that the topology optimization can realize a new BIC for a given {higher} frequency.

We use the same material parameters as previous experiments. Using the topology optimization, we minimize the imaginary part of the resonant wavenumber $\beta=2.10+0.586\ione$ at $\omega=10.0$ of the periodic structure shown in \cref{fig:veri_problem}. To this end, we set the objective functional $J$ as $J=(\mathrm{Im}\,[\beta])^2$ and determine an optimized unit structure within the fixed design domain $[-0.354,0.354]\times [-0.354,0.354]$, so that it exhibits a BIC if $J$ attains the value of zero. The size of the fixed design domain is chosen to avoid violating the well-separated condition (described in \Cref{ss:multi}).

\begin{figure}[h]
  \centering
  \includegraphics[scale=0.27]{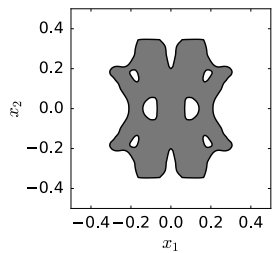}
  \caption{Optimized shape of a unit structure.}
  \label{fig:topt_shape}
 \end{figure}
\begin{figure}[h]
 \centering
 \includegraphics[scale=0.27]{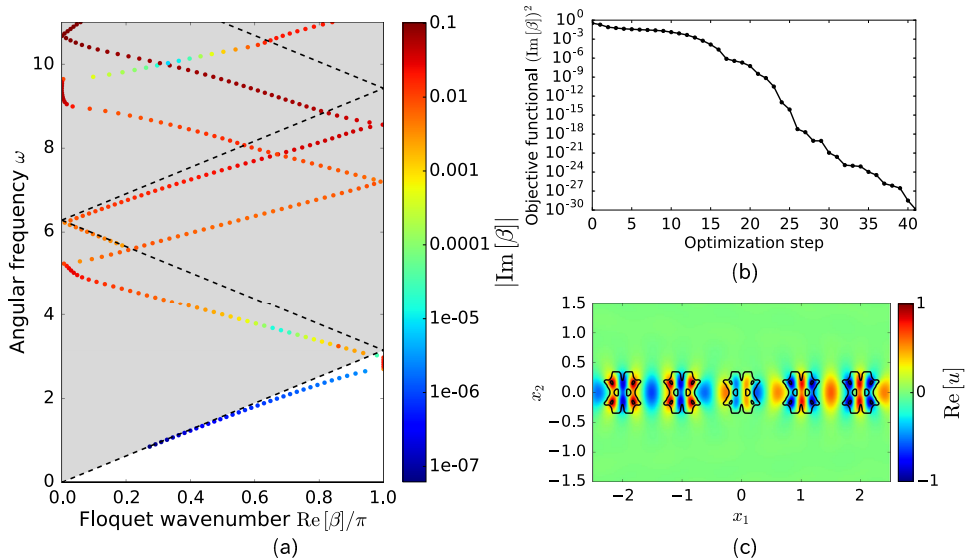}
 \caption{Result of the topology optimization. (a) Band structure for the optimized structure. (b) Convergence history of the optimization.  (c) Mode profile of the BIC at $(\omega,\beta)=(10.0,1.00+0.109\times10^{-14}\ione)$.}
 \label{fig:topt}
\end{figure}
\cref{fig:topt_shape,fig:topt} show the results of the topology optimization. We obtain the optimized shape shown in {\cref{fig:topt_shape}} using the topology optimization for the unit structure. This structure has a resonant wavenumber of $\beta=1.00+0.109\times10^{-14}\ione$ (corresponding to the objective value $J=1.19\times10^{-30}$) at $\omega=10.0$. \cref{fig:topt} (b) shows the convergence history of $J$. The figure shows that the topology optimization successfully decreases the value of $J$. We also conduct a band analysis for the optimized shape and plot the band structure in \cref{fig:topt} (a). From the band structure, we observe that the optimized shape has small imaginary parts around $\omega=10.0$, whereas the initial shape has relatively large imaginary parts (\cref{fig:band_init} (a)). {Although the obtained eigenvalue has a significantly small imaginary part, we cannot guarantee that this is a true BIC because of numerical errors that arise in the BEM and SSM.}

\begin{figure}[h]
\centering
 \includegraphics[scale=0.28]{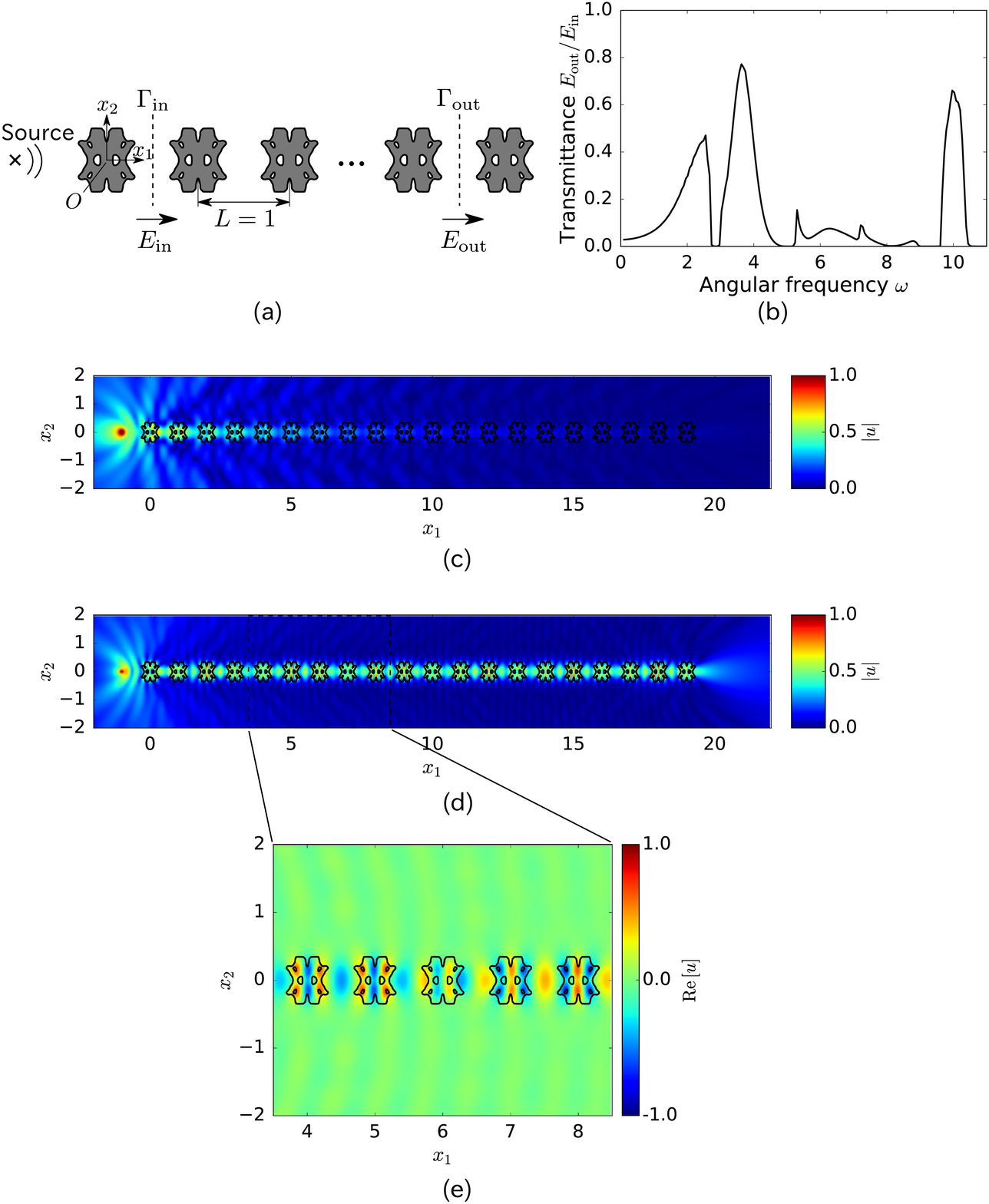}
 \caption{Scattering of a cylindrical wave by the optimized structure. (a) Array of the optimized unit structure; the array comprises 20 scatterers. (b) Transmittance spectrum of guided waves along the structure. (c) Intensity of the total field for $\omega=8.0$. (d) Intensity of the total field for $\omega=10.0$. (e) Real part of the total field for $\omega=10.0$.}
 \label{fig:topt_spectrum}
\end{figure}
To show that the optimized structure supports a guided wave at the desired frequency $\omega=10.0$, we investigate the scattering of the cylindrical wave $H^{(1)}_0(k|\bm{x}-\bm{x}_\mathrm{source}|)$ through the optimized array with source point $\bm{x}_\mathrm{source}=(-1.0,0.0)^T$ (\cref{fig:topt_spectrum}). We compute the energy fluxes $E_\mathrm{in}$ and $E_\mathrm{out}$ across the lines $\Gamma_\mathrm{in}=\{\bm{x}\in\mathbb{R}^2 \mid x_1=0.5,\ -0.5<x_2<0.5\}$ and $\Gamma_\mathrm{out}=\{\bm{x}\in\mathbb{R}^2 \mid x_1=19.5,\ -0.5<x_2<0.5\}$, respectively. \cref{fig:topt_spectrum} (b) shows the plot of the transmittance $E_\mathrm{out}/E_\mathrm{in}$ and frequency. The figure shows that the spectrum has peaks at $\omega=3.6$ and $10.0$, corresponding to the eigenvalues with small imaginary parts in \cref{fig:topt} (a). \Cref{fig:topt_spectrum} (c) and (d) show the total field $u$ when the incident wave illuminates the array for $\omega=8.0$ and $\omega=10.0$, respectively. We also plot the real part of the total field at $\omega=10.0$ in \Cref{fig:topt_spectrum} (e); it shows the similar wave profile to the new BIC (\cref{fig:topt} (c)). These results show that the incident field excites the guided mode (BIC) at $\omega=10.0$ realized by the topology optimization; however, it exhibits no coupling to any guided mode for $\omega=8.0$.

\begin{figure}[h]
\centering
 \includegraphics[scale=0.28]{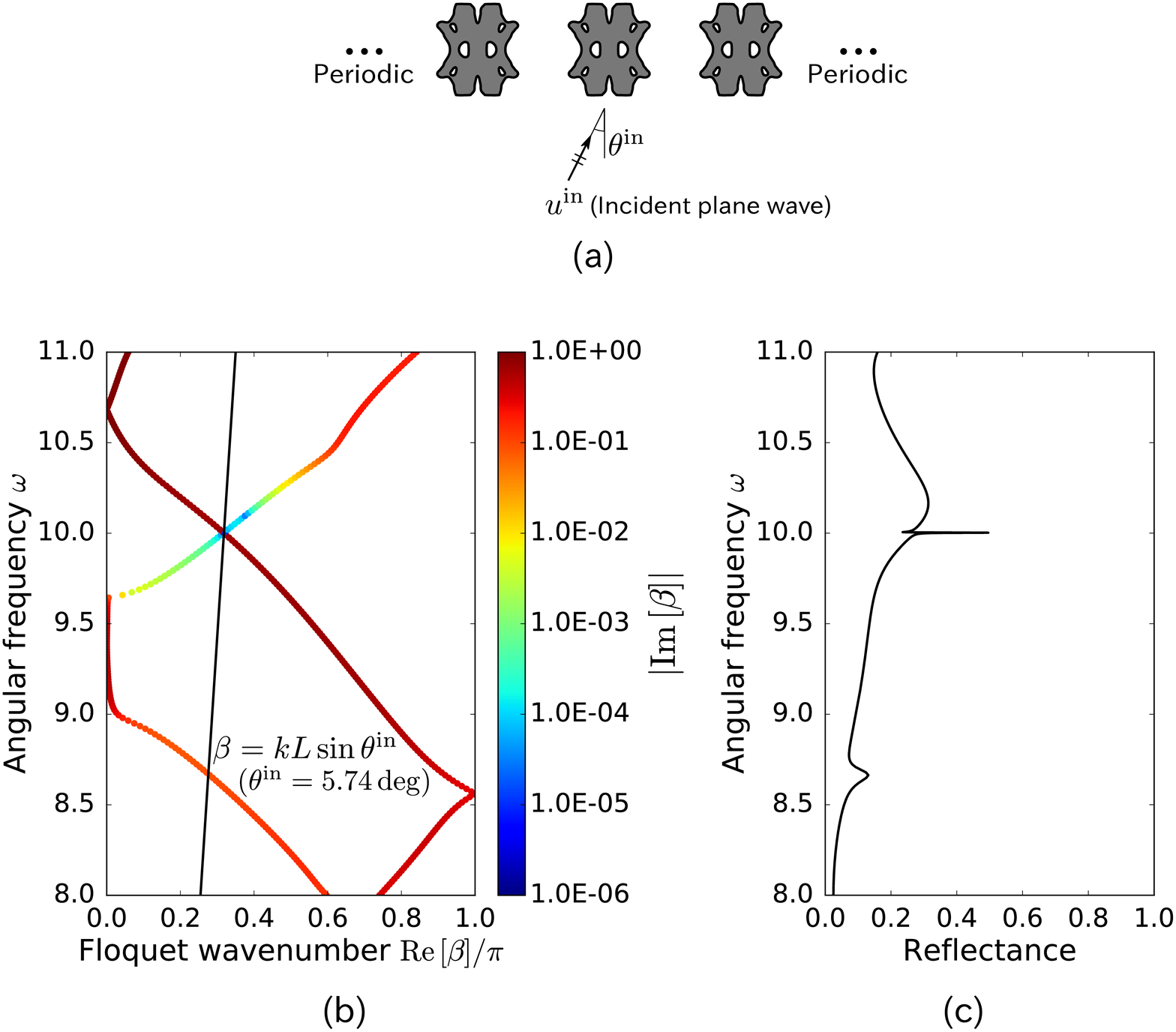}
 \caption{Scattering of a plane wave by the optimized periodic structure. (a) Optimized periodic structure. (b) Band structure. (c) Reflectance spectrum for $\theta^\mathrm{in}=5.74\,\mathrm{deg}$.}
 \label{fig:topt_scattering}
\end{figure}
We can also observe a BIC by finding Wood's anomaly in a scattering analysis \cite{Monticone2017bound}. As shown in \cref{fig:topt_scattering} (a), we analyze the scattering of a plane wave by the optimized structure. The incident angle $\theta^\mathrm{in}$ is given by $\theta^\mathrm{in}=5.74\,\mathrm{deg}$ so that the line $\beta=kL\sin\theta^\mathrm{in}$ crosses the band at $\omega=10.0$ as shown in \cref{fig:topt_scattering} (b). \Cref{fig:topt_scattering} (c) shows the reflectance, which is defined by the downward energy flux divided by the incident energy flux per unit cell, versus the angular frequency. The spectrum exhibits a sharp resonance at $\omega=10.0$, corresponding to the BIC realized through the topology optimization. {Further analyses show that a solution of the scattering problem is not unique at exact BICs \cite{bonnet-bendhia1994guided,shipman2010resonant}.}

\begin{figure}[h]
  \centering
   \includegraphics[scale=0.5]{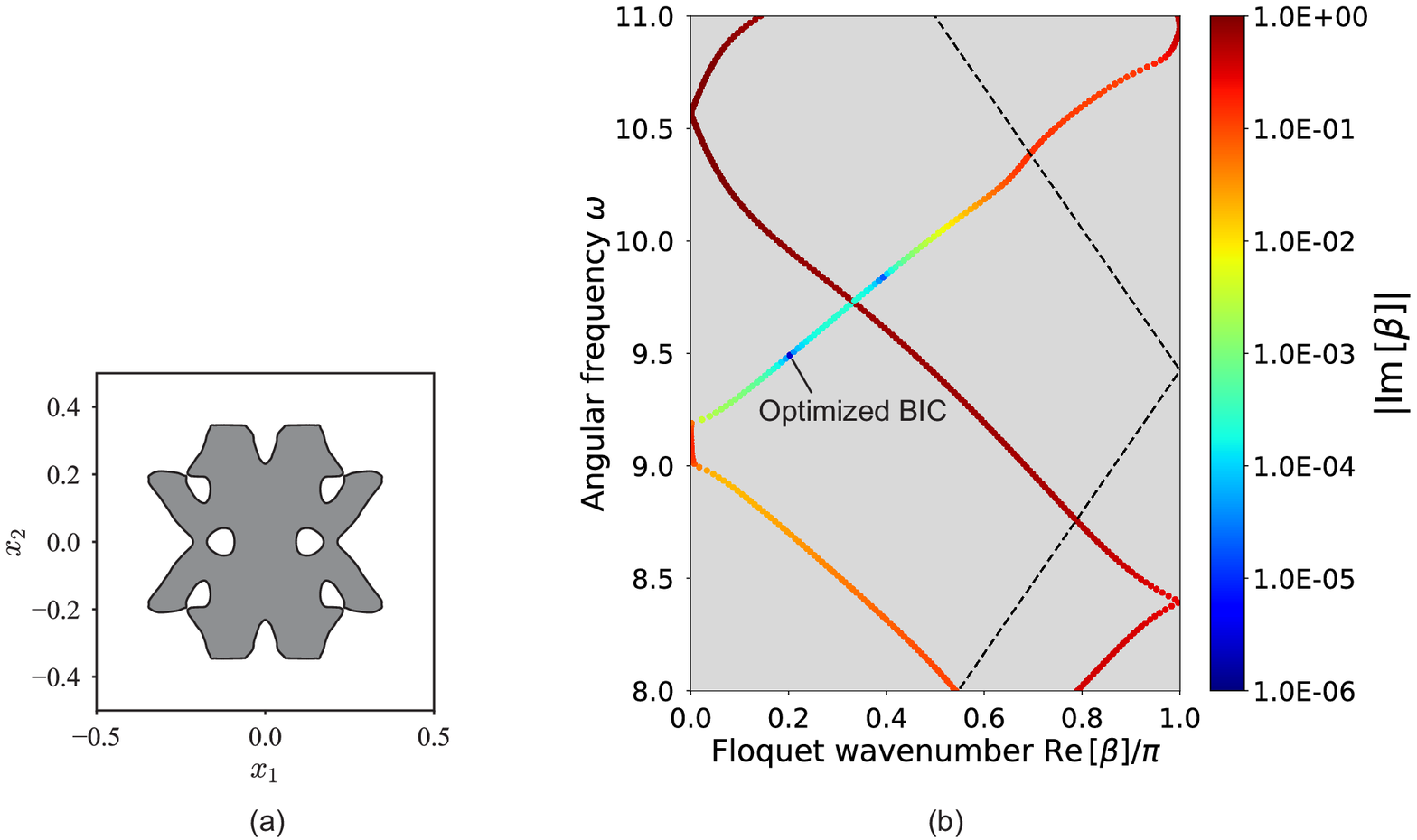}
   \caption{{Optimized structure (a) and band diagram (b) for the target resonant pair $(\omega,\beta)=(9.48,1.15+0.293\ione)$.}}
   \label{fig:revision_result}
  \end{figure}
{In the band diagram shown in \cref{fig:topt} (a), the optimized BIC occurs at the crossing of two bands. To check whether this is a necessary condition for realizing a BIC, we conduct the same topology optimization at the frequency $\omega=9.48$ with resonant wavenumber $\beta=1.15+0.293\ione$ in the initial diagram, shown in \cref{fig:band_init} (a). The optimized geometry and diagram are plotted in \cref{fig:revision_result}. The results show that the obtained BIC with $\beta=0.615+4.22\times10^{-5}\ione$ is not a crossing point in the diagram, meaning that the crossing is not a necessary condition.}

\section{Conclusions}
This study proposed a topology optimization scheme for designing resonant waveguides exhibiting BICs at desired frequencies in the two-dimensional space. We formulated the periodic problem using the scattering matrix of a unit structure and computed resonant wavenumbers using the BEM and SSM. Moreover, we derived a topological derivative of resonant wavenumbers. In the numerical experiments, we first demonstrated that the proposed method determines a resonant wavenumber accurately. Subsequently, we performed a topology optimization to realize a new BIC at a given frequency. {Although we considered Helmholtz' equation for a singly periodic system in two dimensions, the underlying idea, which is the combination of BEM, SSM, and topology optimization, would be applicable to two-dimensional problems and other wave fields governed by Maxwell's equations or elastodynamics.}


\section*{Acknowledgements}
{The authors would like to acknowledge anonymous referees for their valuable comments.} This work was supported by JSPS KAKENHI Grant Numbers JP19J21766 and JP19H00740.


\begin{thebibliography}{10}
  \expandafter\ifx\csname url\endcsname\relax
    \def\url#1{\texttt{#1}}\fi
  \expandafter\ifx\csname urlprefix\endcsname\relax\def\urlprefix{URL }\fi
  \expandafter\ifx\csname href\endcsname\relax
    \def\href#1#2{#2} \def\path#1{#1}\fi
  
  \bibitem{hsu2016bound}
  C.~W. Hsu, B.~Zhen, A.~D. Stone, J.~D. Joannopoulos, M.~Solja{\v c}i{\'c},
    Bound states in the continuum, Nature Reviews Materials 1~(9) (2016) 1--13.
  \newblock \href {http://dx.doi.org/10.1038/natrevmats.2016.48}
    {\path{doi:10.1038/natrevmats.2016.48}}.
  
  \bibitem{vonneumann1929uber}
  J.~{von Neumann}, E.~P. Wigner, \"uber merkw\"urdige diskrete {{Eigenwerte}},
    Physikalische Zeitschrift 30 (1929) 465--467.
  
  \bibitem{plotnik2011experimental}
  Y.~Plotnik, O.~Peleg, F.~Dreisow, M.~Heinrich, S.~Nolte, A.~Szameit, M.~Segev,
    Experimental observation of optical bound states in the continuum, Physical
    Review Letters 107~(18) (2011) 183901.
  \newblock \href {http://dx.doi.org/10.1103/PhysRevLett.107.183901}
    {\path{doi:10.1103/PhysRevLett.107.183901}}.
  
  \bibitem{lee2012observation}
  J.~Lee, B.~Zhen, S.-L. Chua, W.~Qiu, J.~D. Joannopoulos, M.~Solja{\v c}i{\'c},
    O.~Shapira, Observation and differentiation of unique high-{$Q$} optical
    resonances near zero wave vector in macroscopic photonic crystal slabs,
    Physical Review Letters 109~(6) (2012) 067401.
  \newblock \href {http://dx.doi.org/10.1103/PhysRevLett.109.067401}
    {\path{doi:10.1103/PhysRevLett.109.067401}}.
  
  \bibitem{hsu2013observation}
  C.~W. Hsu, B.~Zhen, J.~Lee, S.-L. Chua, S.~G. Johnson, J.~D. Joannopoulos,
    M.~Solja{\v c}i{\'c}, Observation of trapped light within the radiation
    continuum, Nature 499~(7457) (2013) 188--191.
  \newblock \href {http://dx.doi.org/10.1038/nature12289}
    {\path{doi:10.1038/nature12289}}.
  
  \bibitem{linton2007embedded}
  C.~M. Linton, P.~McIver, Embedded trapped modes in water waves and acoustics,
    Wave Motion 45~(1) (2007) 16--29.
  \newblock \href {http://dx.doi.org/10.1016/j.wavemoti.2007.04.009}
    {\path{doi:10.1016/j.wavemoti.2007.04.009}}.
  
  \bibitem{dreisow2009adiabatic}
  F.~Dreisow, A.~Szameit, M.~Heinrich, R.~Keil, S.~Nolte, A.~T{\"u}nnermann,
    S.~Longhi, Adiabatic transfer of light via a continuum in optical waveguides,
    Optics Letters 34~(16) (2009) 2405--2407.
  \newblock \href {http://dx.doi.org/10.1364/OL.34.002405}
    {\path{doi:10.1364/OL.34.002405}}.
  
  \bibitem{yang2014analytical}
  Y.~Yang, C.~Peng, Y.~Liang, Z.~Li, S.~Noda, Analytical perspective for bound
    states in the continuum in photonic crystal slabs, Physical Review Letters
    113~(3) (2014) 037401.
  \newblock \href {http://dx.doi.org/10.1103/PhysRevLett.113.037401}
    {\path{doi:10.1103/PhysRevLett.113.037401}}.
  
  \bibitem{bulgakov2014bloch}
  E.~N. Bulgakov, A.~F. Sadreev, Bloch bound states in the radiation continuum in
    a periodic array of dielectric rods, Physical Review A 90~(5) (2014) 053801.
  \newblock \href {http://dx.doi.org/10.1103/PhysRevA.90.053801}
    {\path{doi:10.1103/PhysRevA.90.053801}}.
  
  \bibitem{kodigala2017lasing}
  A.~Kodigala, T.~Lepetit, Q.~Gu, B.~Bahari, Y.~Fainman, B.~Kant{\'e}, Lasing
    action from photonic bound states in continuum, Nature 541~(7636) (2017)
    196--199.
  \newblock \href {http://dx.doi.org/10.1038/nature20799}
    {\path{doi:10.1038/nature20799}}.
  
  \bibitem{foley2014symmetryprotected}
  J.~M. Foley, S.~M. Young, J.~D. Phillips, Symmetry-protected mode coupling near
    normal incidence for narrow-band transmission filtering in a dielectric
    grating, Physical Review B 89~(16) (2014) 165111.
  \newblock \href {http://dx.doi.org/10.1103/PhysRevB.89.165111}
    {\path{doi:10.1103/PhysRevB.89.165111}}.
  
  \bibitem{yanik2011seeing}
  A.~A. Yanik, A.~E. Cetin, M.~Huang, A.~Artar, S.~H. Mousavi, A.~Khanikaev,
    J.~H. Connor, G.~Shvets, {Hatice Altug}, Seeing protein monolayers with naked
    eye through plasmonic {{Fano}} resonances, Proceedings of the National
    Academy of Sciences 108~(29) (2011) 11784--11789.
  \newblock \href {http://dx.doi.org/10.1073/pnas.1101910108}
    {\path{doi:10.1073/pnas.1101910108}}.
  
  \bibitem{sokolowski1992introduction}
  J.~Sokolowski, J.-P. Zolesio, Introduction to {{Shape Optimization}}: {{Shape
    Sensitivity Analysis}}, Springer {{Series}} in {{Computational Mathematics}},
    {Springer}, {Berlin, Heidelberg}, 1992.
  
  \bibitem{bendsoe2013topology}
  M.~P. Bends{\o}e, O.~Sigmund, Topology Optimization: Theory, Methods, and
    Applications, {Springer Science \& Business Media}, 2013.
  
  \bibitem{dobson1999maximizing}
  D.~C. Dobson, S.~J. Cox, Maximizing band gaps in two-dimensional photonic
    crystals, SIAM Journal on Applied Mathematics 59~(6) (1999) 2108--2120.
  \newblock \href {http://dx.doi.org/10.1137/S0036139998338455}
    {\path{doi:10.1137/S0036139998338455}}.
  
  \bibitem{evans1998trapped}
  D.~Evans, R.~Porter, Trapped modes embedded in the continuous spectrum,
    Quarterly Journal of Mechanics and Applied Mathematics 51~(2) (1998)
    263--274.
  \newblock \href {http://dx.doi.org/10.1093/qjmam/51.2.263}
    {\path{doi:10.1093/qjmam/51.2.263}}.
  
  \bibitem{porter2005embedded}
  R.~Porter, D.~V. Evans, Embedded {{Rayleigh}}\textendash{{Bloch}} surface waves
    along periodic rectangular arrays, Wave Motion 43~(1) (2005) 29--50.
  \newblock \href {http://dx.doi.org/10.1016/j.wavemoti.2005.05.005}
    {\path{doi:10.1016/j.wavemoti.2005.05.005}}.
  
  \bibitem{bennetts2022rayleigh}
  L.~G. Bennetts, M.~A. Peter, Rayleigh\textendash{{Bloch}} waves above the
    cutoff, Journal of Fluid Mechanics 940.
  \newblock \href {http://dx.doi.org/10.1017/jfm.2022.247}
    {\path{doi:10.1017/jfm.2022.247}}.
  
  \bibitem{hu2009understanding}
  J.~Hu, C.~R. Menyuk, Understanding leaky modes: Slab waveguide revisited,
    Advances in Optics and Photonics 1~(1) (2009) 58--106.
  \newblock \href {http://dx.doi.org/10.1364/AOP.1.000058}
    {\path{doi:10.1364/AOP.1.000058}}.
  
  \bibitem{ammari2020perturbation}
  H.~Ammari, A.~Dabrowski, B.~Fitzpatrick, P.~Millien, Perturbation of the
    scattering resonances of an open cavity by small particles. {{Part I}}: The
    transverse magnetic polarization case, Zeitschrift f\"ur angewandte
    Mathematik und Physik 71~(4) (2020) 102.
  \newblock \href {http://dx.doi.org/10.1007/s00033-020-01324-6}
    {\path{doi:10.1007/s00033-020-01324-6}}.
  
  \bibitem{gimbutas2013fast}
  Z.~Gimbutas, L.~Greengard, Fast multi-particle scattering: {{A}} hybrid solver
    for the {{Maxwell}} equations in microstructured materials, Journal of
    Computational Physics 232~(1) (2013) 22--32.
  \newblock \href {http://dx.doi.org/10.1016/j.jcp.2012.01.041}
    {\path{doi:10.1016/j.jcp.2012.01.041}}.
  
  \bibitem{sokolowski1999topological}
  J.~Sokolowski, A.~Zochowski, On the topological derivative in shape
    optimization, SIAM Journal on Control and Optimization 37~(4) (1999)
    1251--1272.
  \newblock \href {http://dx.doi.org/10.1137/S0363012997323230}
    {\path{doi:10.1137/S0363012997323230}}.
  
  \bibitem{isakari2017topology}
  H.~Isakari, T.~Takahashi, T.~Matsumoto, A topology optimisation with level-sets
    of {{B-spline}} surface (in {{Japanese}}), Transactions of the Japan Society
    for Computational Methods in Engineering 17 (2017) 125--130.
  
  \bibitem{burton1971application}
  A.~J. Burton, G.~F. Miller, J.~H. Wilkinson, The application of integral
    equation methods to the numerical solution of some exterior boundary-value
    problems, Proceedings of the Royal Society of London. A. Mathematical and
    Physical Sciences 323~(1553) (1971) 201--210.
  \newblock \href {http://dx.doi.org/10.1098/rspa.1971.0097}
    {\path{doi:10.1098/rspa.1971.0097}}.
  
  \bibitem{zheng2015burton}
  C.-J. Zheng, H.-B. Chen, H.-F. Gao, L.~Du, Is the
    {{Burton}}\textendash{{Miller}} formulation really free of fictitious
    eigenfrequencies?, Engineering Analysis with Boundary Elements 59 (2015)
    43--51.
  \newblock \href {http://dx.doi.org/10.1016/j.enganabound.2015.04.014}
    {\path{doi:10.1016/j.enganabound.2015.04.014}}.
  
  \bibitem{abramowitz1965handbook}
  M.~Abramowitz, I.~A. Stegun, Handbook of Mathematical Functions with Formulas,
    Graphs, and Mathematical Tables, {Dover Publications}, 1965.
  
  \bibitem{martin2006multiple}
  P.~A. Martin, Multiple {{Scattering}}: {{Interaction}} of {{Time-harmonic
    Waves}} with {{N Obstacles}}, {Cambridge University Press}, 2006.
  
  \bibitem{coifman1993fast}
  R.~Coifman, V.~Rokhlin, S.~Wandzura, The fast multipole method for the wave
    equation: A pedestrian prescription, IEEE Antennas and Propagation Magazine
    35~(3) (1993) 7--12.
  \newblock \href {http://dx.doi.org/10.1109/74.250128}
    {\path{doi:10.1109/74.250128}}.
  
  \bibitem{nicorovici1995photonic}
  N.~A. Nicorovici, R.~C. McPhedran, L.~C. Botten, Photonic band gaps for arrays
    of perfectly conducting cylinders, Physical Review E 52~(1) (1995)
    1135--1145.
  \newblock \href {http://dx.doi.org/10.1103/PhysRevE.52.1135}
    {\path{doi:10.1103/PhysRevE.52.1135}}.
  
  \bibitem{linton2006schlomilch}
  C.~M. Linton, Schl\"omilch series that arise in diffraction theory and their
    efficient computation, Journal of Physics A: Mathematical and General 39~(13)
    (2006) 3325--3339.
  \newblock \href {http://dx.doi.org/10.1088/0305-4470/39/13/012}
    {\path{doi:10.1088/0305-4470/39/13/012}}.
  
  \bibitem{porter1999rayleigh}
  R.~Porter, D.~V. Evans, Rayleigh\textendash{{Bloch}} surface waves along
    periodic gratings and their connection with trapped modes in waveguides,
    Journal of Fluid Mechanics 386 (1999) 233--258.
  \newblock \href {http://dx.doi.org/10.1017/S0022112099004425}
    {\path{doi:10.1017/S0022112099004425}}.
  
  \bibitem{otani2008fmm}
  Y.~Otani, N.~Nishimura, An {{FMM}} for periodic boundary value problems for
    cracks for {{Helmholtz}}' equation in {{2D}}, International Journal for
    Numerical Methods in Engineering 73~(3) (2008) 381--406.
  \newblock \href {http://dx.doi.org/10.1002/nme.2077}
    {\path{doi:10.1002/nme.2077}}.
  
  \bibitem{isakari2012calderon}
  H.~Isakari, K.~Niino, H.~Yoshikawa, N.~Nishimura, Calderon's preconditioning
    for periodic fast multipole method for elastodynamics in {{3D}},
    International Journal for Numerical Methods in Engineering 90~(4) (2012)
    484--505.
  \newblock \href {http://dx.doi.org/10.1002/nme.3332}
    {\path{doi:10.1002/nme.3332}}.
  
  \bibitem{nose2014calculation}
  T.~Nose, N.~Nishimura, Calculation of eigenvalues related to 2 dimensional
    periodic boundary value problems for the {{Helmholtz}} equation using the
    {{Sakurai-Sugiura}} method and periodic fast multipole method (in
    {{Japanese}}), Transactions of the Japan Society for Industrial and Applied
    Mathematics 24 (2014) 185--201.
  
  \bibitem{misawa2016fmm}
  R.~Misawa, K.~Niino, N.~Nishimura, An {{FMM}} for waveguide problems of 2-{{D
    Helmholtz}}' equation and its application to eigenvalue problems, Wave Motion
    63 (2016) 1--17.
  \newblock \href {http://dx.doi.org/10.1016/j.wavemoti.2015.12.006}
    {\path{doi:10.1016/j.wavemoti.2015.12.006}}.
  
  \bibitem{ooura1999robust}
  T.~Ooura, M.~Mori, A robust double exponential formula for {{Fourier-type}}
    integrals, Journal of Computational and Applied Mathematics 112~(1) (1999)
    229--241.
  \newblock \href {http://dx.doi.org/10.1016/S0377-0427(99)00223-X}
    {\path{doi:10.1016/S0377-0427(99)00223-X}}.
  
  \bibitem{asakura2009numerical}
  J.~Asakura, T.~Sakurai, H.~Tadano, T.~Ikegami, K.~Kimura, A numerical method
    for nonlinear eigenvalue problems using contour integrals, JSIAM Letters 1
    (2009) 52--55.
  \newblock \href {http://dx.doi.org/10.14495/jsiaml.1.52}
    {\path{doi:10.14495/jsiaml.1.52}}.
  
  \bibitem{nakamoto2017levelsetbased}
  K.~Nakamoto, H.~Isakari, T.~Takahashi, T.~Matsumoto, A level-set-based topology
    optimisation of carpet cloaking devices with the boundary element method,
    Mechanical Engineering Journal 4~(1) (2017) 16--00268.
  \newblock \href {http://dx.doi.org/10.1299/mej.16-00268}
    {\path{doi:10.1299/mej.16-00268}}.
  
  \bibitem{amstutz2006new}
  S.~Amstutz, H.~Andr{\"a}, A new algorithm for topology optimization using a
    level-set method, Journal of Computational Physics 216~(2) (2006) 573--588.
  \newblock \href {http://dx.doi.org/10.1016/j.jcp.2005.12.015}
    {\path{doi:10.1016/j.jcp.2005.12.015}}.
  
  \bibitem{bonnet-bendhia1994guided}
  A.-S. {Bonnet-Bendhia}, F.~Starling, Guided waves by electromagnetic gratings
    and non-uniqueness examples for the diffraction problem, Mathematical Methods
    in the Applied Sciences 17~(5) (1994) 305--338.
  \newblock \href {http://dx.doi.org/10.1002/mma.1670170502}
    {\path{doi:10.1002/mma.1670170502}}.
  
  \bibitem{shipman2010resonant}
  S.~P. Shipman, Resonant scattering by open periodic waveguides, in: Progress in
    {{Computational Physics}}, Vol.~1, {Bentham Science Publishers}, 2010.
  
  \bibitem{yuan2018bound}
  L.~Yuan, Y.~Y. Lu, Bound states in the continuum on periodic structures
    surrounded by strong resonances, Physical Review A 97~(4) (2018) 043828.
  \newblock \href {http://dx.doi.org/10.1103/PhysRevA.97.043828}
    {\path{doi:10.1103/PhysRevA.97.043828}}.
  
  \bibitem{Monticone2017bound}
  F.~Monticone, A.~Al{\`u}, Bound states within the radiation continuum in
    diffraction gratings and the role of leaky modes, New Journal of Physics
    19~(9) (2017) 093011.
  \newblock \href {http://dx.doi.org/10.1088/1367-2630/aa849f}
    {\path{doi:10.1088/1367-2630/aa849f}}.
  
  \end{thebibliography}

\end{document}